\documentclass[english]{amsart}

\usepackage{amsmath}

\usepackage{amssymb}

\usepackage{amscd}
\usepackage{tabularx}

\usepackage{appendix}
\usepackage{pdfsync}

\usepackage[all]{xy}

\newcommand{\F}{\mathbb{F}} 

\DeclareMathOperator{\Frob}{Frob}

 \newcommand{\To}{\longrightarrow}
 \newcommand{\isoto}{\stackrel{\sim}{\To}}
 
   \newcommand{\onto}{\twoheadrightarrow}
     \newcommand{\into}{\hookrightarrow}

 \newcommand{\M}{\mathcal{M}}

 \newcommand{\bigO}{\mathcal{O}}
  \newcommand{\T}{\mathbb{T}}

 \newcommand{\Z}{\mathbb{Z}}
 
 \newcommand{\Q}{\mathbb{Q}}

\newcommand{\mcal}{\mathcal}

 \newcommand{\Gal}{\operatorname{Gal}}
 \newcommand{\GL}{\operatorname{GL}}

 \newcommand{\BrMod}{\operatorname{BrMod}}
 
 \newcommand{\rhobar}{\overline{\rho}}

\newcommand{\End}{\operatorname{End}}
\newcommand{\Hom}{\operatorname{Hom}}

\newcommand{\Aut}{\operatorname{Aut}}

\newcommand{\Fil}{\operatorname{Fil}}

\newcommand{\Sym}{\operatorname{Sym}}
\newcommand{\Mod}{\operatorname{Mod}}

\usepackage{amsthm}

\newcommand{\Qp}{{{\mathbb Q}_p}}
\newcommand{\Fp}{{\mathbb F}_p}

\newcommand{\Qbar}{\overline{\Q}}
\newcommand{\Fpbar}{\overline{\F}_p}

\usepackage[mathscr]{eucal}

\newcommand{\Qpbar}{{\overline{\Q}_p}}

\newtheorem*{thmn}{Theorem}

 \newtheorem{thm}{Theorem}[subsection]
 
 \newtheorem{lemma}[thm]{Lemma}
 \newtheorem{prop}[thm]{Proposition}
 \theoremstyle{definition}
 \newtheorem{defn}[thm]{Definition}
\newtheorem{rem}[thm]{Remark}
 \theoremstyle{remark}
 
 \numberwithin{equation}{subsection}

\theoremstyle{definition}

\usepackage{hyperref}
\begin{document}

\title{On the weights of mod $p$ Hilbert modular forms}
\author{Toby Gee}


\begin{abstract}We prove many cases of a conjecture of Buzzard, Diamond and
Jarvis on the possible weights of mod $p$ Hilbert modular forms, by making
use of modularity lifting theorems and computations in $p$-adic Hodge
theory.\end{abstract}
\maketitle
\tableofcontents
\section{Introduction}If a representation$$\overline{\rho}:G_\Q\to\GL_2(\overline{\F}_p)$$ is
continuous, odd, and irreducible, then a conjecture of Serre (now a theorem of Khare-Wintenberger and Kisin) predicts that
$\overline{\rho}$ is modular. More precisely, Serre predicted a minimal
weight $k(\overline{\rho})$ and a minimal level $N(\overline{\rho})$ for a
modular form giving rise to $\overline{\rho}$. 

It is natural to try to extend these results to totally real fields $F$. The
natural generalisation of Serre's conjecture is to conjecture that if
$$\overline{\rho}:G_F\to\GL_2(\overline{\F}_p)$$is continuous, irreducible
and totally odd, then it is modular (in the sense that it arises from a
Hilbert modular form). It is straightforward to generalise the definition of
$N(\overline{\rho})$ to this setting, and there has been much progress on
``level-lowering'' for Hilbert modular forms. It is, however, much harder to
generalise the definition of $k(\overline{\rho})$. For example, there is no
longer a total ordering on the weights, and the $p$-adic Hodge theory is
much more complicated than in the classical case.

Suppose that $p$ is unramified in $F$. Recently (see \cite{bdj}), Buzzard,
Diamond and Jarvis have proposed a conjectural set $W(\overline{\rho})$ of
weights attached to $\overline{\rho}$, from which in the classical case one can deduce the weight part of Serre's conjecture (see \cite{bdj} for more details). In this
paper we prove many cases of a closely related conjecture (we work with a definite, rather than indefinite quaternion algebra; as we discuss below, it should be straightforward to prove the corresponding results in the setting of \cite{bdj}). To be precise, a weight is
an irreducible $\Fpbar$-representation of $\GL_2(\bigO_F/p)$, and such a
representation factors as a tensor product
$$\otimes_{v|p}\sigma_{\vec{a},\vec{b}}$$where $\vec{a}$, $\vec{b}$ are
$[k_v:\F_p]$-tuples indexed by embeddings
$\tau:k_v\hookrightarrow\overline{\F}_p$, and $0\leq a_\tau\leq p-1$,
$1\leq b_\tau\leq p$. Then we say that a weight is \emph{regular} if
in fact $2\leq b_\tau\leq p-2$ for all $\tau$. Our main theorem
requires a technical condition which we prefer to define later, that of
a weight being partially ordinary of type $I$ for $\rhobar$, $I$ a set
of places of $F$ dividing $p$; see section \ref{2}.

\begin{thmn}Suppose that $\overline{\rho}$ is modular, that $p\ge 5$, and
  that $\overline{\rho}|_{G_{F(\zeta_{p})}}$ is irreducible. Then if
  $\sigma$ is a regular weight and $\overline{\rho}$ is modular of
  weight $\sigma$ then $\sigma\in W(\overline{\rho})$. Conversely, if
  $\sigma\in W(\overline{\rho})$ and $\sigma$ is non-ordinary for
  $\rhobar$, then $\overline{\rho}$ is modular of weight $\sigma$. If
  $\sigma$ is partially ordinary of type $I$ for $\rhobar$ and
  $\overline{\rho}$ has a partially ordinary modular lift of type $I$
  then $\overline{\rho}$ is modular of weight $\sigma$.\end{thmn}

Before we discuss the proof, we make some remarks about the
assumptions in the theorem. The assumption that $\overline{\rho}$ is
modular is essential to our methods. The assumption that $p\ge 5$ is
needed in order for there to be any regular weights at all; it is
possible that this could be relaxed in future work, as there is no
essential obstruction to the application of the techniques that we
employ in characteristic $3$. In characteristic $2$ the results from
$2$-adic Hodge theory that we would require have not yet been
developed in sufficient generality, but this too does not appear to be
an insurmountable difficulty. The assumption
that $\overline{\rho}|_{G_{F(\zeta_p)}}$ is irreducible, and the
assumption on partial ordinarity, are needed in order to apply $R=T$
theorems.

The main idea of our proof is the same as that for our proof of a companion
forms theorem for totally real fields (see \cite{gee051}), namely that we
use a  lifting theorem to construct lifts
of $\overline{\rho}$ satisfying certain local properties at places $v|p$,
and then use a modularity lifting theorem of Kisin to prove that these
representations are modular. In fact, Kisin's theorem is not general enough
for our applications, and we need to use the main theorem of \cite{gee052}.
The arguments are much more complicated than those in \cite{gee051} because
we need to construct liftings with more delicate local properties; rather
than just considering ordinary lifts, we must consider potentially
Barsotti-Tate lifts of specified type. 

The other complication which intervenes is that the connection between being
modular of a certain weight and having a lift of a certain type is rather
subtle, and this is the reason for our hypothesis that the weight be
regular. One needs to consider many liftings for each weight, and we have
only obtained the necessary combinatorial results in the case where the
weight is regular. However, while these results appear to hold for most
non-regular weights, there are cases where they do not hold, so it seems
that it is not possible to give a general proof that the list of weights is
correct by simply considering the types of potentially Barsotti-Tate lifts. It is possible to give a complete proof in the case where $p$ splits completely in $F$, and we do this in \cite{gee061}.

We now outline the structure of the paper. Rather than working with
the ``geometric'' conventions of \cite{bdj}, we prefer to work with
more ``arithmetic'' ones. In particular, we work with automorphic
forms on definite quaternion algebras. We set out our conventions in
section \ref{2}, and we state the appropriate reformulation of the
conjectures of \cite{bdj} here. In section \ref{reducible} we carry
out the required local analysis in the case where the local
representation is reducible. Sections \ref{brmod} and \ref{deform} use
Breuil modules and strongly divisible modules to determine when
reducible representations arise as the generic fibres of certain
finite flat group schemes. In section \ref{fl} we relate these finite
flat group schemes to certain crystalline representations considered
in \cite{bdj}, and in section \ref{reducibletypes} we prove the
necessary combinatorial results relating types and regular weights.

We then repeat this analysis in the  irreducible case in
section \ref{irreducible}, and finally in section \ref{global} we combine
these results with the lifting theorems mentioned above to deduce our main
results. Firstly, we use our local results to show that if $\rhobar$ is modular of weight $\sigma$ with $\sigma$ regular, then $\sigma\in W(\rhobar)$. For each regular weight $\sigma\in W(\rhobar)$ we then produce a modular lift of $\rhobar$ which is potentially Barsotti-Tate of a specific type, so that $\rhobar$ must be modular of some weight occurring in the mod $p$ reduction of this type. We then check that $\sigma$ is the only element of $W(\rhobar)$ occurring in this reduction, so that $\rhobar$ is modular of weight $\sigma$, as required. In fact, we do not quite do this; the combinatorics is slightly more involved, and we are forced to make use of a notion of a ``weakly regular'' weight. See section \ref{global} for the details.

It is a pleasure to thank Fred Diamond for numerous helpful
discussions regarding this work; without his patient advice this paper
could never have been written. We would like to thank David Savitt for
pointing out several errors and omissions in an earlier version of
this paper, and for writing \cite{sav06}. We would like to thank
Florian Herzig for pointing out an inconsistency between our
conventions and those of \cite{bdj}, which led to the writing of
section \ref{2}. We are extremely grateful to Xavier Caruso and
Christophe Breuil for their many helpful comments and corrections; in
particular, the material in section \ref{fl} owes a considerable debt
to Caruso's efforts to correct a number of inaccuracies, and the proof
of Lemma \ref{nicebasis} is based on an argument of his. We would also
like to thank the anonymous referee for a careful reading, and for
pointing out a number of serious errors in an earlier version of the
paper.

\section{Definitions}\label{2}\subsection{}Rather than use the
conventions of \cite{bdj}, we choose to state a closely related
variant of their conjectures by working on totally definite quaternion
algebras. This formulation is more suited to applications to
modularity lifting theorems, and indeed to the application of
modularity lifting theorems to proving cases of the conjecture.

We begin by recalling some standard facts from the theory of
quaternionic modular forms; see either \cite{taymero}, section 3 of
\cite{kis04} or section 2 of \cite{kis06} for more details, and in
particular the proofs of the results claimed below. We will follow
Kisin's approach closely. We fix throughout this paper an algebraic
closure $\Qbar$ of $\Q$, and regard all algebraic extensions of $\Q$
as subfields of $\Qbar$. For each prime $p$ we fix an algebraic
closure $\Qpbar$ of $\Qp$, and we fix an embedding
$\Qbar\into\Qpbar$. In this way, if $v$ is a finite place of a number
field $F$, we have a homomorphism $G_{F_v}\into G_F$.

 Let $F$ be a totally real field in which
$p>2$ is unramified, and let $D$ be a quaternion algebra with center
$F$ which is ramified at all infinite places of $F$ and at a set
$\Sigma$ of finite places, which contains no places above $p$. Fix a
maximal order $\bigO_D$ of $D$ and for each finite place
$v\notin\Sigma$ fix an isomorphism $(\bigO_D)_v\isoto
M_2(\bigO_{F_v})$. For any finite place $v$ let $\pi_v$ denote a
uniformiser of $F_v$.

Let $U=\prod_v U_v\subset (D\otimes_F\mathbb{A}^f_F)^\times$ be a compact subgroup, with each $U_v\subset (\bigO_D)^\times_v$. Furthermore, assume that $U_v=(\bigO_D)_v^\times$ for all $v\in\Sigma$, and that $U_v=\GL_2(\bigO_{F_v})$ if $v|p$.

Take $A$ a topological $\Z_p$-algebra. For each $v|p$, fix a
continuous representation $\sigma_v:U_v\to\Aut(W_{\sigma_v})$ with
$W_{\sigma_v}$ a finite free $A$-module. Write
$W_\sigma=\otimes_{v|p,A}W_{\sigma_v}$ and let
$\sigma=\prod_{v|p}\sigma_v$. We regard $\sigma$ as a representation
of $U$ in the obvious way (that is, we let $U_v$ act trivially if
$v\nmid p$). Fix also a character
$\psi:F^\times\backslash(\mathbb{A}^f_F)^\times\to A^\times$ such that for any
place $v$ of $F$, $\sigma|_{U_v\cap\bigO_{F_v}^\times}$ is
multiplication by $\psi^{-1}$. Then we can think of $W_\sigma$ as a
$U(\mathbb{A}_F^f)^\times$-module by letting $(\mathbb{A}_F^f)^\times$
act via $\psi^{-1}$.

Let $S_{\sigma,\psi}(U,A)$ denote the set of continuous
functions $$f:D^\times\backslash(D\otimes_F\mathbb{A}_F^f)^\times\to
W_\sigma$$ such that for all $g\in (D\otimes_F\mathbb{A}_F^f)^\times$
we have $$f(gu)=\sigma(u)^{-1}f(g)\text{ for all }u\in
U,$$ $$f(gz)=\psi(z)f(g)\text{ for all
}z\in(\mathbb{A}_F^f)^\times.$$We can write
$(D\otimes_F\mathbb{A}_F^f)^\times=\coprod_{i\in I}D^\times
t_iU(\mathbb{A}_F^f)^\times$ for some finite index set $I$ and some
$t_i\in(D\otimes_F\mathbb{A}_F^f)^\times$. Then we
have $$S_{\sigma,\psi}(U,A)\isoto\oplus_{i\in
  I}W_\sigma^{(U(\mathbb{A}_F^f)^\times\cap t_i^{-1}D^\times
  t_i)/F^\times},$$the isomorphism being given by the direct sum of
the maps $f\mapsto\{f(t_i)\}$. From now on we make the following
assumption:$$\text{For all
}t\in(D\otimes_F\mathbb{A}_F^f)^\times\text{ the group
}(U(\mathbb{A}_F^f)^\times\cap t^{-1}D^\times t)/F^\times=1.$$ One can
always replace $U$ by a subgroup (obeying the assumptions above) for
which this holds (cf. section 3.1.1 of \cite{kis07}). Under this
assumption, which we make from now on, $S_{\sigma,\psi}(U,A)$ is a
finite free $A$-module, and the functor $W_\sigma\mapsto
S_{\sigma,\psi}(U,A)$ is exact in $W_\sigma$.

We now define some Hecke algebras. Let $S$ be a set of finite places containing $\Sigma$, the places dividing $p$, and the primes of $F$ such that $U_v$ is not a maximal compact subgroup of $ D_v^\times$. Let $\T^{\operatorname{univ}}_{S,A}=A[T_v]_{v\notin S}$ be the commutative polynomial ring in the formal variables $T_v$. Consider the left action of $(D\otimes_F\mathbb{A}_F^f)^\times$ on $W_\sigma$-valued functions on $(D\otimes_F\mathbb{A}_F^f)^\times$ given by $(gf)(z)=f(zg)$. For each finite place $v$ of $F$ we fix a uniformiser $\pi_v$ of $F_v$. Then we make $S_{\sigma,\psi}(U,A)$ a $\T^{\operatorname{univ}}_{S,A}$-module by letting $T_v$ act via the double coset  $U\bigl(
\begin{smallmatrix}
	\pi_v&0\\0&1
\end{smallmatrix}
\bigr)U$. These are independent of the choices of $\pi_v$. We will write $\T_{\sigma,\psi}(U,A)$ or $\T_{\sigma,\psi}(U)$ for the image of $\T^{\operatorname{univ}}_{S,A}$ in $\End S_{\sigma,\psi}(U,A)$.

Let $\mathfrak{m}$ be a maximal ideal of $\T^{\operatorname{univ}}_{S,A}$. We say that $\mathfrak{m}$ is in the support of $(\sigma,\psi)$ if $S_{\sigma,\psi}(U,A)_\mathfrak{m}\neq 0$. Now let $\bigO$ be the ring of integers in  $\Qpbar$, with residue field $\F=\Fpbar$, and suppose that $A=\bigO$ in the above discussion, and that $\sigma$ has open kernel. Consider a maximal ideal $\mathfrak{m}\subset\T^{\operatorname{univ}}_{S,\bigO}$ which is induced by a maximal ideal of $\T_{\sigma,\psi}(U,\bigO)$. Then there is a semisimple Galois representation $\overline{\rho}_\mathfrak{m}:G_F\to\GL_2(\F)$ associated to $\mathfrak{m}$ which is characterised up to equivalence by the property that if $v\notin S$ and $\Frob_v$ is an arithmetic Frobenius at $v$, then the trace of $\overline{\rho}_\mathfrak{m}(\Frob_v)$ is the image of $T_v$ in $\F$. 

We are now in a position to define what it means for a Galois
representation to be modular of some weight. Let $v|p$ be a place of
$F$, let $F_v$ have ring of
integers $\bigO_v$ and residue field $k_v$, and let $\sigma$ be an
irreducible $\overline{\F}_p$-representation of
$G:=\prod_{v|p}\GL_2(k_v)$. We also denote by $\sigma$ the
representation of $\prod_{v|p}\GL_2(\bigO_{v})$ induced by the
surjections $\bigO_{v}\onto k_v$.
\begin{defn}
  We say that an irreducible representation
  $\overline{\rho}:G_F\to\GL_2(\overline{\F}_p)$ is modular of weight
  $\sigma$ if for some $D$, $S$, $U$, $\psi$, and $\mathfrak{m}$ as
  above we have $S_{\sigma,\psi}(U,\F)_\mathfrak{m}\neq 0$ and
  $\overline{\rho}_\mathfrak{m}\cong\overline{\rho}$.
\end{defn}

We now show how one can gain information about the weights associated to a particular Galois representation by considering lifts to characteristic zero.
\begin{lemma}
  \label{432} Let
  $\psi:F^{\times}\backslash(\mathbb{A}_{F}^f)^{\times}\to\bigO^{\times}$
  be a continuous character, and write $\overline{\psi}$ for the
  composite of $\psi$ with the projection $\bigO^{\times}\to
  \F^{\times}$. Fix a representation $\sigma$ of $\prod_{v|p}U_v$ on a finite free
  $\bigO$-module $W_{\sigma}$, and an irreducible representation
  $\sigma'$ on a finite free $\F$-module $W_{\sigma'}$. Suppose that
  for each $v|p$ we have
  $\sigma|_{U_v\cap\bigO_{F_v}^\times}=\psi^{-1}|_{U_v\cap\bigO_{F_v}^\times}$
  and
  $\sigma'|_{U_v\cap\bigO_{F_v}^\times}=\overline{\psi}^{-1}|_{U_v\cap\bigO_{F_v}^\times}$
	
	Let $\mathfrak{m}$ be a maximal ideal of $\mathbb{T}_{S,\bigO}^{\textrm{univ}}$.
	
	Suppose that $W_{\sigma'}$ occurs as a $\prod_{v|p}U_v$-module subquotient of ${W}_{\overline{\sigma}}:=W_\sigma\otimes\F$. If $\mathfrak{m}$ is in the support of $(\sigma',\overline{\psi})$, then $\mathfrak{m}$ is in the support of $(\sigma,\psi)$.
	
	Conversely, if $\mathfrak{m}$ is in the support of $(\sigma,\psi)$, then $\mathfrak{m}$ is in the support of $(\sigma',\overline{\psi})$ for some irreducible $\prod_{v|p}U_v$-module subquotient $W_{\sigma'}$ of ${W}_{\overline{\sigma}}$.
\end{lemma}
\begin{proof}
	The first part is proved just as in Lemma 3.1.4 of \cite{kis04}, and the second part follows from Proposition 1.2.3 of \cite{as86}.
\end{proof}

We note a special case of this result, relating the existence of potentially Barsotti-Tate lifts of a particular tame type to information about Serre weights. Firstly, we recall some particular representations of $\GL_{2}(k_{v})$. 
For any pair of distinct characters
$\chi_1$, $\chi_2:k_{v}^\times\to\bigO^\times$ we let
$I(\chi_1,\chi_2)$ denote the irreducible $(q+1)$-dimensional $\Qpbar$-representation of
$\GL_2(k_{v})$ induced from the
character of $B$ (the upper triangular matrices in $\GL_2(k_{v})$) given by $$\left(\begin{array}{cc}
  x & w \\
  0 & y \\
\end{array}
\right)\mapsto\chi_1(x)\chi_2(y).$$ We let
$\sigma_{\chi_{1},\chi_{2}}$ denote the representation of
$\GL_{2}(k_{v})$ on an $\bigO$-lattice in $I(\chi_1,\chi_2)$; we also
regard this as a representation of $\GL_{2}(\bigO_{v})$ via the
natural projection. Let $\tau(\sigma_{\chi_{1},\chi_{2}})$ be the
inertial type $\chi_{1}\oplus\chi_{2}$ (regarded as a representation
of $I_{F_{v}}$ via local class field theory, normalised so that a
uniformiser corresponds to a geometric Frobenius element).

Let $k_{v}'$ be the quadratic extension of $k_{v}$.  For any character
$\theta:k_{v}'^\times\to\bigO^{\times}$ which does not factor through
the norm $k_{v}'^\times\to k_v^\times$, there is an irreducible
$(q-1)$-dimensional cuspidal representation $\Theta(\theta)$ of
$\GL_2(k_v)$ (see Section 1 of \cite{dia05} for the definition of
$\Theta(\theta)$).  Let $\sigma_{\Theta(\theta)}$ denote the
representation of $\GL_{2}(k_{v})$ on an $\bigO$-lattice in
$\Theta(\theta)$; we also regard this as a representation of
$\GL_{2}(\bigO_{v})$ via the natural projection. Let $q_{v}$ be the
cardinality of $k_{v}$, and let $\tau(\sigma_{\Theta(\theta)})$ be the
inertial type $\theta\oplus\theta^{q_{v}}$ (again regarded as a
representation of $I_{F_{v}}$ via local class field theory).

\begin{defn}\label{defn:barsotti tate lifts} Let $\tau$ be an inertial
  type, and let $v|p$ be a place of $F$. We say that a lift $\rho$ of
  $\rhobar|_{G_{F_v}}$ is \emph{potentially Barsotti-Tate of type} $\tau$ if $\rho$ is potentially
  Barsotti-Tate, has determinant a finite order
  character of order prime to $p$ times the cyclotomic character, and
  the corresponding Weil-Deligne representation (see Appendix B of \cite{cdt}), when restricted to
  $I_{F_v}$, is isomorphic to $\tau$. \end{defn}

\begin{lemma}\label{lem:local langlands version of lifting}For each
  $v|p$, fix a representation $\sigma_{v}$ of the type just considered
  (that is, isomorphic to $\sigma_{\chi_{1},\chi_{2}}$ or to
  $\sigma_{\Theta(\theta)}$). Let $\tau_{v}=\tau(\sigma_{v})$ be the
  corresponding inertial type. Suppose that $\rhobar$ is modular of
  weight $\sigma$, and that $\sigma$ is a
  $\prod_{v|p}\GL_{2}(k_{v})$-subquotient of
  $\otimes_{v|p}\sigma_{v}\otimes_{\bigO}\F$. Then $\rhobar$ lifts to
  a modular Galois representation which is potentially Barsotti-Tate
  of type $\tau_{v}$ for each $v|p$.
	
	Conversely, suppose that $\rhobar$ lifts to a modular Galois representation which is potentially Barsotti-Tate of type $\tau_{v}$ for each $v|p$. Then  $\rhobar$ is modular of weight $\sigma$ for some  $\prod_{v|p}\GL_{2}(k_{v})$-subquotient $\sigma$ of $\otimes_{v|p}\sigma_{v}\otimes_{\bigO}\F$.
	
\end{lemma}
\begin{proof}This follows from Lemma \ref{432}, the Jacquet-Langlands correspondence, and the compatibility of the local and global Langlands correspondences at places dividing $p$ (see \cite{kis06}).
	
\end{proof}
We now state a conjecture on Serre weights, following \cite{bdj}. Note
that our conjecture is only valid for regular weights (a notion which
we will define shortly); there are some additional complications when
dealing with non-regular weights. Let $\rhobar:G_F\to\GL_2(\Fpbar)$ be
modular. We propose a conjectural set of regular weights $W(\rhobar)$
for $\rhobar$. 

In fact, for each place $v|p$ we propose a set of weights $W(\rhobar|_{G_{F_v}})$, and we define $$W(\rhobar):=\bigl\{\otimes_{v|p}\sigma_v | \sigma_v\in W(\rhobar|_{G_{F_v}})\bigr\}.$$

Let $S_v$ be the set of embeddings $k_v\into \Fpbar$.  A weight for $\GL_2(k_v)$ is
an isomorphism class of irreducible $\overline{\F}_p$-representations of
$\GL_2(k_{v})$, which automatically contains one of the form
$$\sigma_{\vec{a},\vec{b}}=\otimes_{\tau\in
S_{v}}\det{}^{a_\tau}\Sym^{b_\tau-1}k_{v}^2\otimes_\tau\overline{\F}_p,$$with $0\leq
a_\tau\leq p-1$ and $1\leq b_\tau\leq p$ for each $\tau\in S_{v}$. We demand
further that some $a_\tau<p-1$, in which case the representations
$\sigma_{\vec{a},\vec{b}}$ are pairwise non-isomorphic.

\begin{defn}We say that a weight $\sigma_{\vec{a},\vec{b}}$ is \emph{regular} if
$2\leq b_\tau\leq p-2$ for all $\tau$. 	We say that it is \emph{weakly regular} if
	$1\leq b_\tau\leq p-1$ for all $\tau$.\end{defn}

For each $\tau\in S_v$ we have the fundamental character $\omega_\tau$ of $I_{F_v}$ given by composing $\tau$ with the homomorphism $I_{F_v}\to k_v^\times$ given by local class field theory, normalised so that uniformisers correspond to geometric Frobenius elements. Let $k'_v$ denote the quadratic extension of $k_v$. Let $S'_{v}$ denote the set of embeddings $\sigma:k_v'\into\Fpbar$, and let $\omega_\sigma$ denote the fundamental character corresponding to $\sigma$.

Suppose firstly that $\rhobar|_{G_{F_v}}$ is irreducible.  There is a natural $2-1$ map $\pi:S'_v\to S_{v}$ given by restriction to $k_{v}$, and we say that a subset $J\subset S'_{v}$ is a \emph{full subset} if $|J|=|\pi(J)|=|S_{v}|$. Then we have

\begin{defn}Let $\sigma_{\vec{a},\vec{b}}$ be a regular weight for $\GL_2(k_v)$. Then $\sigma_{\vec{a},\vec{b}}\in W(\rhobar|_{G_{F_v}})$ if and only if there exists a full subset $J\subset S'_v$ such that $$\rhobar|_{I_{F_v}}\sim\prod_{\tau\in S_v}\omega_\tau^{a_\tau}\left(\begin{array}{cc}
  \prod_{\sigma\in J}\omega_{\sigma}^{b_{\sigma|_{k_{v}}}}   & 0 \\
  0 & \prod_{\sigma\notin J}\omega_{\sigma}^{b_{\sigma|_{k_{v}}}} \\
\end{array}\right).$$\end{defn}
Suppose now that $\rhobar|_{G_{F_v}}$ is reducible, say $\rhobar|_{G_{F_v}}\sim\bigl(\begin{smallmatrix}\psi_1&*\\0&\psi_2\end{smallmatrix}\bigr)$. We define the set $W(\rhobar|_{G_{F_v}})$ in two stages. Firstly, define a set  $W(\rhobar|_{G_{F_v}})'$ of regular weights as follows.

\begin{defn}\label{defn:weight set in reducible case}Let $\sigma_{\vec{a},\vec{b}}$ be a regular weight for
  $\GL_2(k_v)$. Then 
  $\sigma_{\vec{a},\vec{b}}\in W(\rhobar|_{G_{F_v}})'$ if and only if there exists $J\subset S_v$
  such that $\psi_1|_{I_{F_v}}=\prod_{\tau\in
    S_v}\omega_\tau^{a_\tau}\prod_{\tau\in J}\omega_\tau^{b_\tau}$ and
  $\psi_{2}|_{I_{F_v}}=\prod_{\tau\in
    S_v}\omega_\tau^{a_\tau}\prod_{\tau\notin
    J}\omega_\tau^{b_\tau}$. We say that $\sigma_{\vec{a},\vec{b}}\in
  W(\rhobar|_{G_{F_v}})'$ is \emph{ordinary for} $\rhobar$ if
  furthermore $J=S_v$ or $J=\emptyset$ (note that the set $J$ is
  uniquely determined, because $\sigma_{\vec{a},\vec{b}}$ is regular).\end{defn}

Suppose that we have a regular weight $\sigma_{\vec{a},\vec{b}}\in
W(\rhobar|_{G_{F_v}})'$ and a corresponding subset $J\subset S_v$. We
now define crystalline lifts $\widetilde{\psi}_1$,
$\widetilde{\psi}_2$ of $\psi_1$, $\psi_2$. If $\psi:G_{F_v}\to\Qpbar^\times$ is a crystalline
character, and $\tau:F_{v}\into \Qpbar$, we say that the
Hodge-Tate weight of $\psi$ with respect to $\tau$ is the $i$ for
which
$gr^{-i}((\psi\otimes_{\Qp}B_{dR})^{G_{F_v}}\otimes_{\Qpbar\otimes_\Qp
  F_v,1\otimes\tau}\Qpbar)\neq 0$. Then we demand that for some fixed
Frobenius element $\Frob_v$ of $G_{F_{v}}$,
$\widetilde{\psi}_i(\Frob_v)$ is the Teichm\"{u}ller lift of
$\psi_i(\Frob_v)$, and that:
\begin{itemize}\item $\widetilde{\psi}_1$ is crystalline, and the Hodge-Tate weight of $\widetilde{\psi}_1$ with respect to $\tau$ is $a_\tau+b_\tau$ if $\tau\in J$, and $a_\tau$ if $\tau\notin J$.
\item $\widetilde{\psi}_2$ is crystalline, and the Hodge-Tate weight of $\widetilde{\psi}_2$ with respect to $\tau$ is $a_\tau+b_\tau$ if $\tau\notin J$, and $a_\tau$ if $\tau\in J$.\end{itemize}
The existence and uniqueness (for our fixed choice of $\Frob_v$) of $\widetilde{\psi}_1$, $\widetilde{\psi}_2$ is straightforward (see \cite{bdj}). Then we have

\begin{defn}$\sigma_{\vec{a},\vec{b}}\in W(\rhobar|_{G_{F_v}})$ if and
  only if $\rhobar|_{G_{F_v}}$ has a lift to a crystalline
  representation
  $\bigl(\begin{smallmatrix}\widetilde{\psi}_1&*\\0&\widetilde{\psi}_2\end{smallmatrix}\bigr)$.\end{defn}Note
that by remark 3.10 of \cite{bdj}, and the regularity of
$\sigma_{\vec{a},\vec{b}}$, this definition is independent of the
choice of $\Frob_v$.

For future reference, we say that a weight $\sigma$ is partially
ordinary of type $I$ for $\rhobar$ if $I$ is the set of places $v|p$
for which $\sigma_v$ is ordinary for $\rhobar$. We say that $\rhobar$
has a partially ordinary modular lift of type $I$ if it has a
potentially Barsotti-Tate modular lift which is potentially ordinary
at precisely the places in $I$.

\subsection{Relation to the Buzzard-Diamond-Jarvis conjecture}Our
conjectured sets of regular weights are exactly the same as the
regular weights predicted in \cite{bdj}. However, they work with
indefinite quaternion algebras rather than the definite ones of this
paper, and in the absence of a mod $p$ Jacquet-Langlands
correspondence our results do not automatically prove cases of their
conjectures. That said, our arguments are for the most part purely
local, with the only global input being in characteristic zero, where
one does have a Jacquet-Langlands corresponence. In particular, given
the analogue of Lemma \ref{lem:local langlands version of lifting} in
the setting of \cite{bdj} (cf. Proposition 2.10 of \cite{bdj}) our
arguments will go over unchanged to their setting.

\section{Local analysis - the reducible case}\label{reducible}
\subsection{Breuil Modules}\label{brmod}Let $p>2$ be prime, let $k$ be a
finite extension of $\F_p$, let $K_0=W(k)[1/p]$, and let $K$ be a
finite Galois totally tamely ramified extension of $K_0$, of degree
$e$. Fix a subfield $M$ of $K_0$, and assume that there is a
uniformiser $\pi$ of $\bigO_{K}$ such that $\pi^{e}\in M$, and fix
such a $\pi$.  Since $K/M$ is tamely ramified (and automatically
Galois), the category of Breuil modules with coefficients and descent
data is easy to describe (see \cite{sav06}). Let $k\in[2,p-1]$ be an
integer (there will never be any ambiguity in our two uses of the
symbol $k$, one being a finite field and the other a positive integer). Let $E$ be a finite extension field of $\Fp$. The category
$\BrMod^{k-1}_{dd,M}$ consists of quintuples
$(\mathcal{M},\mathcal{M}_{k-1},\phi_{k-1},\hat{g},N)$ where:

\begin{itemize}\item $\mathcal{M}$ is a finitely generated
  $(k\otimes_{\F_p}E)[u]/u^{ep}$-module, free over $k[u]/u^{ep}$.
\item $\M_{k-1}$ is a $(k\otimes_{\F_p}E)[u]/u^{ep}$-submodule of $\M$
  containing $u^{e(k-1)}\M$.
\item $\phi_{k-1}:\M_{k-1}\to\M$ is $E$-linear and $\phi$-semilinear
  (where $\phi:k[u]/u^{ep}\to k[u]/u^{ep}$ is the $p$-th power map)
  with image generating $\M$ as a
  $(k\otimes_{\F_p}E)[u]/u^{ep}$-module.
\item $N:\M\to u\M$ is $(k\otimes_{\F_{p}}E)$-linear and satisfies
  $N(ux)=uN(x)-ux$ for all $x\in\M$,
  $u^{e}N(\M_{k-1})\subset\M_{k-1}$, and
  $\phi_{k-1}(u^{e}N(x))=(-\pi^e/p)N(\phi_{k-1}(x))$ for all $x\in\M_{k-1}$.
\item $\hat{g}:\M\to\M$ are additive bijections for each
  $g\in\Gal(K/M)$, preserving $\M_{k-1}$, commuting with the $\phi_{k-1}$-,
  $E$-, and $N$-actions, and satisfying $\hat{g}_1\circ
  \hat{g}_2=\widehat{g_1\circ g}_2$ for all $g_1,g_2\in\Gal(K/M)$, and
  $\hat{1}$ is the identity. Furthermore, if $a\in
  k\otimes_{\F_{p}}E$, $m\in\M$ then
  $\hat{g}(au^{i}m)=g(a)((g(\pi)/\pi)^{i}\otimes 1)u^{i}\hat{g}(m)$.\end{itemize}

We will omit $M$ from the notation in the case $M=K_{0}$. We write
$\BrMod_{dd,M}=\BrMod^{1}_{dd,M}$. The category $\BrMod_{dd,M}$ is
equivalent to the category of finite flat group schemes over
$\mathcal{O}_K$ together with an $E$-action and descent data on the
generic fibre from $K$ to $M$ (this equivalence depends on $\pi$). In
this case it follows from the other axioms that there is always a
unique $N$ which satisfies the required properties, and we will frequently omit
the details of this operator when we are working in the case $k=2$. In
section \ref{fl} we will also use the case $k=p-1$, and here we will
make the operators $N$ explicit.

We choose in this paper (except in section \ref{fl}) to adopt the
conventions of \cite{bm} and \cite{sav04}, rather than those of
\cite{bcdt}; thus rather than working with the usual contravariant
equivalence of categories, we work with a covariant version of it, so
that our formulae for generic fibres will differ by duality and a
twist from those following the conventions of \cite{bcdt}.  To be
precise, we obtain the associated $G_{M}$-representation (which we
will refer to as the generic fibre) of an object of $\BrMod_{dd}$ via
the functor $T_{st,2}^{M}$, which is defined in section 4 of
\cite{sav04}.

Let $\rho:G_{K_0}\to\GL_2(E)$ be a continuous representation. We
assume from now on that $E$ contains $k$.  Suppose for the rest of
this section that $\rho$ is reducible but not scalar, say $\rho\sim
\bigl(\begin{smallmatrix}\psi_1&*\\0&\psi_2\end{smallmatrix}\bigr)$. Fix
$\pi=(-p)^{1/(p^r-1)}$, where $r=[k:\F_p]$, and fix $K=K_0(\pi)$, so
that $\pi$ is a uniformiser of $\mathcal{O}_K$, the ring of integers
of $K$. By class field theory $\psi_{1}|_{I_K}$ and
$\psi_{2}|_{I_{K}}$ are trivial.

We fix some general notation
for elements of $\BrMod_{dd}$. Let $S$ denote the set of embeddings $\tau:k\hookrightarrow E$. We have an isomorphism
$k\otimes_{\F_p}E\isoto\oplus_{S}E_{\tau}$, where $E_{\tau}:=k\otimes_{k,\tau}E$, and we let
$\epsilon_\tau$ denote the idempotent corresponding to the embedding
$\tau$. Then any element $\mathcal{M}$ of $\BrMod_{dd}$
can be decomposed into $E[u]/u^{ep}$-modules
$\mathcal{M}^{\tau}:=\epsilon_\tau\M$, $\tau\in S$, so that
$\hat{g}:\M^\tau\to\M^\tau$, and
$\phi_1:\M^\tau_1\to\M^{\tau\circ\phi^{-1}}$, so that $\mathcal{M}$ is
free over $(k\otimes_{\Fp} E)[u]/u^{ep}$. We now write
$S=\{\tau_1,\dots,\tau_r\}$, numbered so that
$\tau_{i+1}=\tau_{i}\circ\phi^{{-1}}$, where we identify $\tau_{r+1}$
with $\tau_1$. In fact, it will often be useful to consider the indexing set of $S$ to be $\Z/r\Z$, and we will do so without further comment.

 Fix $J\subset S$. We wish to single out
particular representations $\rho$ depending on $J$. Firstly, we need some
notation. Recall that (as in appendix
B of \cite{cdt}) if $\rho':G_{K_0}\to\GL_2(\mathcal{O}_L)$ is potentially
Barsotti-Tate, where $L$ is a finite extension of $W(E)[1/p]$, then there is
a Weil-Deligne representation
$WD(\rho'):W_{K_0}\to\GL_2(\overline{\Q}_p)$ associated to $\rho'$,
and we say that $\rho'$ has type $WD(\rho')|_{I_{K_0}}$.

\begin{defn}We say that $\rho$ has a \emph{lift of type} $J$ if there is a
representation $\rho':G_{K_0}\to\GL_2(\mathcal{O}_L)$ lifting $\rho$, where
$L$ is a finite extension of $W(E)[1/p]$, such that $\rho'$ becomes
Barsotti-Tate over $K$, with $\varepsilon^{-1}\det\rho'$ equal to the Teichm\"{u}ller
lift of $\varepsilon^{-1}\det\rho$ (with $\varepsilon$ denoting the
cyclotomic character) and $\rho'$ has type
$\widetilde{\psi}_1|_{I_{K_{0}}}\prod_{\tau\in
J}\widetilde{\omega}_\tau^{-p}\oplus\widetilde{\psi}_2|_{I_{K_{0}}}\prod_{\tau\notin
J}\widetilde{\omega}_\tau^{-p}$. Here a tilde denotes the Teichm\"{u}ller
lift.\end{defn}

\begin{defn}For any subset $H\subset S$, we say that an element $\M$ of $\BrMod_{dd}$ is
of class $H$ if it is of rank one, and for all $\tau\in S$ we can choose a basis $e_\tau$ of
$\M^{\tau}$ such that $\M_1^{\tau}$ is generated by $u^{j_\tau}e_\tau$,
where
$$j_\tau=\left\{\begin{array}{c}
                 0\text{ if }\tau\circ\phi^{-1}\notin H \\
                 e\text{ if }\tau\circ\phi^{-1}\in H \\
               \end{array}\right.$$
	
\end{defn}

\begin{defn}We say that an element $\M$ of $\BrMod_{dd}$ is of type $J$ if $\M$ is an
extension of an element of class $J^c$ by an element of class $J$, and we
say that $\rho$ has a model of type $J$ if there is an element of
$\BrMod_{dd}$ of type $J$ with generic fibre $\rho$.
	
\end{defn}
 We will also refer to finite flat group schemes with descent data as being of class $J$ or of type $J$ if they correspond to Breuil modules with descent data of this kind. The notions of having a
model of type $J$ and having a lift of type $J$ are closely related,
although not in general equivalent. We will see in section \ref{subsec: models of type J} that in sufficiently generic cases,
if $\rho$ has a model of type $J$ then it has a lift of type $J$, and in
section \ref{reducibletypes} we prove a partial converse (see Proposition \ref{prop:liftimpliesmodel}).

\subsection{Strongly divisible modules}\label{deform}In this section we prove  that if $\rho$ has a model of type $J$
then it has a lift of type $J$. We begin by recalling the definition
and basic properties of strongly divisible modules from \cite{sav04}. For
the purpose of giving these definitions we return briefly to the general
setting of $K_0$ an unramified finite extension of $\Q_p$ and $K$ a totally
tamely ramified Galois extension of $K_0$ of degree $e$, with uniformiser $\pi$, satisfying $\pi^{e}\in M$ for some subfield $M$ of $K_{0}$.

Let $L$ be a finite extension of $\Q_p$ with ring of integers $\bigO_L$ and
residue field $E$. Let $S_{K}$ be the ring
$$\left\{\sum_{j=0}^\infty r_j\frac{u^j}{\lfloor j/e\rfloor !}\text{, }r_j\in
  W(k)\text{, }r_j\to 0\ p\text{-adically as }j\to\infty\right\},$$
and let $S_{K,\bigO_{L}}=S_{K}\otimes_{\Z_{p}}\bigO_{L}$. Let
$\Fil^{1}S_{K,\bigO_L}$ be the $p$-adic completion of the ideal
generated by $E(u)^j/j!$, $j\geq 1$, where $E(u)$ is the minimal
polynomial of $\pi$ over $K_0$. Let $\phi:S_{K,\bigO_L}\to
S_{K,\bigO_L}$ be the unique $\bigO$-linear, $W(k)$-semilinear ring
homomorphism with $\phi(u)=u^p$, and let $N$ be the unique
$W(k)\otimes \bigO_L$-linear derivation such that $N(u)=-u$ (so that
$N\phi=p\phi N$). One can check that
$\phi(\Fil^{1}S_{K,\bigO_L})\subset p S_{K,\bigO_L}$, and we define
$\phi_{1}:\Fil^{1}S_{K,\bigO_L}\to S_{K,\bigO_L}$ by
$\phi_{1}=(\phi|_{\Fil^{1}S_{K,\bigO_L}})/p$. One can check (see
section 4 of \cite{sav04}) that if $I$ is an ideal of $\bigO_{L}$,
then
$IS_{K,\bigO_L}\cap\Fil^{1}S_{K,\bigO_L}=I\Fil^{1}S_{K,\bigO_L}$. We
give $S_{K}$ an action of $\Gal(K/M)$ via ring isomorphisms via the
usual action on $W(k)$, and by letting $\hat{g}(u)=(g(\pi)/\pi)u$. We
extend this action $\bigO_L$-linearly to $S_{K,\bigO_L}$.

We now define the category $\bigO_{L}-\Mod^1_{cris,dd,M}$, the category of strongly divisible
$\bigO_{L}$-modules with descent data from $K$ to $M$.

\begin{defn}\label{strongdivis}A  strongly divisible
$\bigO_{L}$-module with descent data from $K$ to $M$ is a finitely generated free
$S_{K,\bigO_L}$-module $\M$, together with a sub-$S_{K,\bigO_L}$-module $\Fil^{1}\M$
and a map $\phi:\M\to\M$, and additive bijections $\hat{g}:\M\to\M$ for each
$g\in\Gal(K/M)$, satisfying the following
conditions:\begin{enumerate}\item $\Fil^{1}\M$ contains
$(\Fil^{1}S_{K,\bigO_L})\M$, \item $\Fil^{1}\M\cap I\M=I\Fil^{1}\M$ for all
ideals $I$ in $\bigO_{L}$,\item $\phi(sx)=\phi(s)\phi(x)$ for $s\in S_{K,\bigO_L}$ and
$x\in\M$,\item $\phi(\Fil^{1}\M)$ is contained in $p\M$ and
generates it over $S_{K,\bigO_L}$,\item $\hat{g}(sx)=\hat{g}(s)\hat{g}(x)$ for all $s\in S_{K,\bigO_L}$,
$x\in \M$, $g\in\Gal(K/M)$,\item
$\hat{g}_1\circ\hat{g}_2=\widehat{g_1\circ g}_2$ for all $g_1$,
$g_2\in\Gal(K/M)$,\item $\hat{g}(\Fil^{1}\M)\subset\Fil^{1}\M$ for all
$g\in\Gal(K/M)$, and \item $\phi$ commutes with $\hat{g}$ for all
$g\in\Gal(K/M)$.\end{enumerate}

\end{defn}

Note that it is not immediately obvious that this definition is equivalent to Definition 4.1 of \cite{sav04}, as we have made no mention of the operator $N$ of \emph{loc. cit.} However, since $\bigO_{L}$ is finite over $\Z_{p}$, it follows from part (1) of Proposition 5.1.3 of \cite{MR1804530} that any such operator $N$ is unique. The existence of an operator $N$ satisfying all of the conditions of Definition 4.1 of \cite{sav04} except possibly for $\bigO_{L}$-linearity follows from the argument at the beginning of section 3.5 of \cite{sav04}. To check $\bigO_{L}$-linearity it is enough (by $\Z_{p}$-linearity) to check that $N$ is compatible with the action of the units in $\bigO_{L}$, but this is clear from the uniqueness of $N$.

By Proposition 4.13 of \cite{sav04} (and the remarks immediately preceding it), there is a functor $T^{M}_{st,2}$ from the category
$\bigO_L-\Mod^1_{cris,dd,M}$ to the category of $G_{M}$-stable $\bigO_{L}$-lattices in representations of $G_{M}$ which become Barsotti-Tate on restriction to $G_{K}$. This functor preserves dimensions in the obvious sense.

Recall also from section 4.1 of \cite{sav04})
that there is a functor $T_0$, compatible with
$T_{st,2}^{M}$, from $\bigO_L-\Mod^1_{cris,dd,M}$ to $\BrMod_{dd,M}$. The functor $T_0$ is given by
$\M\mapsto(\M/\mathfrak{m}_{L}\M)\otimes_{S_K} k[u]/u^{ep}$.

\subsection{Models of type $J$}\label{subsec: models of type J}
We now wish to discuss the relationships between models of type $J$ and lifts of type $J$. With an eye to our future applications, we will often make a simplifying assumption.

\begin{defn}
	Say that $\rho$ is \emph{$J$-regular} if $\psi_{1}\psi_{2}^{-1}|_{I_{K_{0}}}=\prod_{\tau\in J}\omega_{\tau}^{b_{\tau}}\prod_{\tau\in J^{c}}\omega_{\tau}^{-b_{\tau}}$ for some $2\leq b_{\tau}\leq p-2$.
\end{defn}

Suppose now that $\rho$ has a model of type $J$. Recall that this
means that, with the notation of Section \ref{brmod}, we can write
down a Breuil module $\M$ with descent data whose generic fibre is
$\rho$, which is an extension of a Breuil module with descent data
$\mathcal{B}$ by a Breuil module with descent data $\mathcal{A}$,
where $\mathcal{A}$ is of class $J$ and $\mathcal{B}$ is of class
$J^{c}$. Let $\psi'_{i}$ denote $\psi_{i}|_{I_{K_{0}}}$ regarded as a
character of $\Gal(K/K_{0})$. By Theorem 3.5 and Example 3.7 of
\cite{sav06} we see that we can choose bases for $\mathcal{A}$ and
$\mathcal{B}$ so that they take the following form:
$$\mathcal{A}^{\tau_i}=E[u]/u^{ep}\cdot e_{\tau_i}$$
$$\mathcal{A}_1^{\tau_i}=E[u]/u^{ep}\cdot u^{j_{\tau_i}} e_{\tau_i}$$
$$\phi_1(u^{j_{\tau_i}} e_{\tau_i})=(a^{-1})_{i}e_{\tau_{i+1}} $$ 
$$\hat{g}(e_{\tau_i})=\left(\left(\psi_1'\prod_{\sigma\in
J}\omega_\sigma^{-p}\right)(g)\right)e_{\tau_i}$$
$$\mathcal{B}^{\tau_i}=E[u]/u^{ep}\cdot \overline{f}_{\tau_i}$$
$$\mathcal{B}_1^{\tau_i}=E[u]/u^{ep}\cdot u^{e-j_{\tau_i}} \overline{f}_{\tau_i}$$
$$\phi_1(u^{e-j_{\tau_i}} \overline{f}_{\tau_i})=(b^{-1})_{i}\overline{f}_{\tau_{i+1}} $$ 
$$\hat{g}(\overline{f}_{\tau_i})=\left(\left(\psi_2'\prod_{\sigma\notin
J}\omega_\sigma^{-p}\right)(g)\right)\overline{f}_{\tau_i}$$
where $a$, $b\in E^{\times}$, the notation $(x)_{i}$ means $x$ if $i=1$ and $1$ otherwise, and
$$j_{\tau_i}=\left\{\begin{array}{c}
                 e\text{ if }\tau_{i+1}\in J \\
                 0\text{ if }\tau_{i+1}\notin J. \\
               \end{array}\right. $$
We now seek to choose a basis for $\M$ extending the basis $\{e_{\tau}\}$ for $\mathcal{A}$. Such a basis will be given by lifting the $\overline{f}_{\tau}$ to elements $f_{\tau}$ (where we mean lifting under the map $e_{\tau}\mapsto 0$).
\begin{lemma}\label{nicebasis}Assume that $\rho$ is $J$-regular and
  has a model $\M$ of type $J$. Then for some choice of basis, we can write $$\M^{\tau_i}=E[u]/u^{ep}\cdot e_{\tau_i}+E[u]/u^{ep}\cdot f_{\tau_i}$$
$$\M_1^{\tau_i}=E[u]/u^{ep}\cdot u^{j_{\tau_i}} e_{\tau_i}+E[u]/u^{ep}\cdot (u^{e-j_{\tau_i}}f_{\tau_i}+\lambda_{\tau_i}u^{i_{\tau_i}}e_{\tau_i})$$
$$\phi_1(u^{j_{\tau_i}} e_{\tau_i})=(a^{-1})_ie_{\tau_{i+1}} $$ $$\phi_1(u^{e-j_{\tau_i}}f_{\tau_i}+\lambda_{\tau_i}u^{i_{\tau_i}}e_{\tau_i})=(b^{-1})_if_{\tau_{i+1}}$$
$$\hat{g}(e_{\tau_i})=\left(\left(\psi_1'\prod_{\sigma\in
J}\omega_\sigma^{-p}\right)(g)\right)e_{\tau_i}$$
$$\hat{g}(f_{\tau_i})=\left(\left(\psi_2'\prod_{\sigma\notin
J}\omega_\sigma^{-p}\right)(g)\right)f_{\tau_i}$$where $\lambda_{\tau_i}\in
E$, with $\lambda_{\tau_i}=0$ if $\tau_{i+1}\notin J$, the $i_{\tau_i}$ are such
that $\M_1$ is Galois-stable and $0\leq i_{\tau_i}\leq e-1$, and
$$j_{\tau_i}=\left\{\begin{array}{c}
                 e\text{ if }\tau_{i+1}\in J \\
                 0\text{ if }\tau_{i+1}\notin J. \\
               \end{array}\right. $$ \end{lemma}
\begin{proof}Assume firstly that $J\neq S$, and choose $k$ so that $\tau_{k+1}\notin J$. One can lift $\overline{f}_{\tau_{k}}$ to an element $f_{\tau_{k}}$ of $\phi_{1}(\M^{\tau_{k-1}})$, and in fact one can choose $f_{\tau_{k}}$ so that for all $g\in \Gal(K/{K_{0}})$ we have 			$$\hat{g}(f_{\tau_k})=\left(\left(\psi_2'\prod_{\sigma\notin
			J}\omega_\sigma^{-p}\right)(g)\right)f_{\tau_k}$$(the obstruction to doing this is easily checked to vanish, as the degree of $K/K_{0}$ is prime to $p$). As $\tau_{k+1}\notin J$, we have $j_{\tau_{k}}=0$, so that $e_{\tau_{k}}$ and $u^{e}f_{\tau_{k}}$ must generate $\M^{\tau_{k}}_{1}$.
			
Now, suppose inductively that for some $i$ we have chosen $f_{\tau_{i}}$ and $\lambda_{\tau_{i}}$ so that $\M_1^{\tau_i}$ is generated by  $u^{j_{\tau_i}} e_{\tau_i}$ and $(u^{e-j_{\tau_i}}f_{\tau_i}+\lambda_{\tau_i}u^{i_{\tau_i}}e_{\tau_i})$. Then we put $f_{\tau_{i+1}}=\phi_{1}(u^{e-j_{\tau_i}}f_{\tau_i}+\lambda_{\tau_i}u^{i_{\tau_i}}e_{\tau_i})/(b^{-1})_{i}$. Then $f_{\tau_{i+1}}$ is a lift of $\overline{f}_{\tau_{i+1}}$, and the commutativity of $\phi_{1}$ and the action of $\Gal(K/K_{0})$ ensures that $$\hat{g}(f_{\tau_{i+1}})=\left(\left(\psi_2'\prod_{\sigma\notin
						J}\omega_\sigma^{-p}\right)(g)\right)f_{\tau_{i+1}}.$$ Then the fact that $\M_{1}$ is $\Gal(K/K_{0})$-stable ensures that for some $\lambda_{t_{i+1}}\in E$ we must have that $u^{j_{\tau_{i+1}}} e_{\tau_{i+1}}$ and $(u^{e-j_{\tau_{i+1}}}f_{\tau_{i+1}}+\lambda_{\tau_{i+1}}u^{i_{\tau_{i+1}}}e_{\tau_{i+1}})$ generate $\M_{1}^{\tau_{i+1}}$, and of course if $\tau_{i+2}\notin J$ we can take $\lambda_{\tau_{i+1}}=0$.
						
So, beginning at $k$ we inductively define $f_{\tau_{i}}$ and $\lambda_{\tau_{i}}$ for all $i$, which automatically satisfy all the required properties, except that we do not know that $$\phi_1(u^{e-j_{\tau_{k-1}}}f_{\tau_{k-1}}+\lambda_{\tau_{k-1}}u^{i_{\tau_{k-1}}}e_{\tau_{k-1}})=(b^{-1})_{k-1} f_{\tau_{k}}.$$However, because $k+1\notin J$, we may replace $f_{\tau_{k}}$ with $\phi_1(u^{e-j_{\tau_{k-1}}}f_{\tau_{k-1}}+\lambda_{\tau_{k-1}}u^{i_{\tau_{k-1}}}e_{\tau_{k-1}})/(b^{-1})_{k-1}$ without altering the fact that $$\phi_1(u^{e}f_{\tau_{k}})=(b^{-1})_{k}f_{\tau_{k+1}},$$ so we are done.
						
Suppose now that $J=S$. Then we may carry out a similar inductive
procedure starting with $\tau_{1}$, and we again define $f_{\tau_{i}}$
and $\lambda_{\tau_{i}}$ for all $i$, satisfying all the required
properties, except that we do not know
that $$\phi_{1}(f_{\tau_{r}}+\lambda_{\tau_{r}}u^{i_{\tau_{r}}}e_{\tau_{r}})=f_{\tau_{1}}.$$
We wish to redefine $f_{\tau_{1}}$ to be
$\phi_{1}(f_{\tau_{r}}+\lambda_{\tau_{r}}e_{\tau_{r}})$, and we claim
that doing so does not affect the
relation $$\phi_{1}(f_{\tau_{1}}+\lambda_{\tau_{1}}u^{i_{\tau_{1}}}e_{\tau_{1}})=b^{-1}f_{\tau_{2}}.$$
To see this, note that we are modifying $f_{\tau_{1}}$ by a multiple
of $e_{\tau_{1}}$ which is in the image of $\phi_{1}$, which by
considering the action of $\Gal(K/K_{0})$ must in fact be of the form
$\theta u^{pi_{\tau_{r}}}e_{\tau_{1}}$, with $\theta\in E$ and
$pi_{\tau_{r}}\equiv i_{\tau_{1}}\textrm{ mod }e.$ Now, the assumption
that $\rho$ is $S$-regular means that
$i_{\tau_{1}}=e-\sum_{l=1}^{r}p^{r-l}(b_{\tau_{l+1}}-1)\equiv
-b_{\tau_{1}}\textrm{ mod }p$, with $2\leq b_{\tau_{l}}\leq p-2$. Now,
if we write $pi_{\tau_{r}}=i_{\tau_{1}}+me$, we see that $m\equiv
i_{\tau_{1}}\equiv -b_{\tau_{1}}\textrm{ mod }p$, and since $2\leq
b_{\tau_{1}}\leq p-2$ we see that $m\geq 2$. But then $\phi_{1}(\theta
u^{pi_{\tau_{r}}}e_{\tau_{1}})=\phi_{1}(\theta
u^{i_{\tau_{1}}+(m-1)e}u^{e}e_{\tau_{1}})$ is divisible by
$u^{p(m-1)e}$ and is thus $0$, as required.

\end{proof}

\begin{thm}\label{localdef}Assume that $\rho$ is $J$-regular and has a model of type $J$. Then $\rho$ has a lift of type $J$, which is potentially ordinary if and only if $J=S$ or $J=\emptyset$. 
\end{thm}
\begin{proof}We will write down an element $\M_{J}$ of
  $W(E)-\Mod_{cris,dd,K_{0}}$ such that $T_{0}(\M_{J})=\M$, where $\M$
  is as in Lemma \ref{nicebasis}. We can write $S_{K,W(E)}$
  as $\oplus_{\tau\in S} S_{K}$, and we then define
	
	$$\M_J^{\tau_i}=S_{K,W(E)}\cdot e_{\tau_i}+S_{K,W(E)}\cdot f_{\tau_i}$$
	 $$\hat{g}(e_{\tau_i})=\left(\left(\widetilde{\psi}_1'\prod_{\sigma\in
	J}\widetilde{\omega}_\sigma^{-p}\right)(g)\right)e_{\tau_i}$$
	$$\hat{g}(f_{\tau_i})=\left(\left(\widetilde{\psi}_2'\prod_{\sigma\notin
	J}\widetilde{\omega}_\sigma^{-p}\right)(g)\right)f_{\tau_i}$$
		If  $\tau_{i+1}\in J$,
	$$\Fil^1\M_J^{\tau_i}=\Fil^1S_{K,W(E)}\cdot\M^{\tau_i}_J+S_{K,W(E)}\cdot (f_{\tau_i}+\tilde{\lambda}_{\tau_i}u^{i_{\tau_i}}e_{\tau_i})$$ $$\phi(e_{\tau_i})=(\tilde{a}^{-1})_i e_{\tau_{i+1}}$$ $$\phi(f_{\tau_i}+\tilde{\lambda}_{\tau_i}u^{i_{\tau_i}}e_{\tau_i})=(\tilde{b}^{-1})_i pf_{\tau_{i+1}}$$
		If  $\tau_{i+1}\notin J$,
	$$\Fil^1\M_J^{\tau_i}=\Fil^1S_{K,W(E)}\cdot\M^{\tau_i}_J+S_{K,W(E)}\cdot e_{\tau_i}$$ $$\phi(e_{\tau_i})=(\tilde{a}^{-1})_i p e_{\tau_{i+1}}$$ $$\phi(f_{\tau_i})=(\tilde{b}^{-1})_i f_{\tau_{i+1}}$$
	
	 Here a tilde denotes a Teichm\"{u}ller lift. 
	
         Firstly we verify that this really is an element of
         $W(E)-\Mod^1_{cris,dd,K_{0}}$. Of the properties in
         Definition \ref{strongdivis}, the only non-obvious points are
         that $\Fil^1\M_J\cap I\M_J=I\Fil^1\M_J$ for all ideals $I$ of
         $\bigO_{L}$, and that $\phi(\Fil^1\M_J)$ is contained in
         $p\M_J$ and generates it over $S_{K,W(E)}$. But these are
         both straightforward; that $\Fil^1\M_J\cap I\M_J=I\Fil^1\M_J$
         follows at once from the definition of $\Fil^1\M_J$ and the
         corresponding assertion for $S_{K}$, and that
         $\phi(\Fil^1\M_J)$ is contained in $p\M_J$ and generates it
         over $S_{K,W(E)}$ follows by inspection and the corresponding
         assertions for $S_K$.
	
         It is immediate from the definition of $T_0$ that
         $T_0(\M_J)\simeq\M$. To see that $T^{K_0}_{st,2}(\M_J)$ is a
         lift of $\rho$ of type $J$, note firstly that the Hodge-Tate
         weights of $T^{K_0}_{st,2}(\M_J)$ can be read off from the
         form of the filtration, exactly as in the last two paragraphs
         of the proof of Theorem 6.1 of \cite{geesavquatalg}. This
         shows that the determinant is a finite order character times
         the cyclotomic character, and it also shows that the
         representation is potentially ordinary if and only if $J=S$
         or $J=\emptyset$. That the lift is of type $J$ is then
         immediate from the form of the $\Gal(K/K_{0})$-action and
         Proposition 5.1 of \cite{geesavquatalg}.\end{proof}

\subsection{Breuil modules and Fontaine-Laffaille theory}\label{fl}In this section we
relate the notion of having a model of type $J$ to that of possessing a certain crystalline lift. Suppose as usual that $\rho\sim
\bigl(\begin{smallmatrix}\psi_1&*\\0&\psi_2\end{smallmatrix}\bigr)$, and
that we can write $\psi_1|_{I_{K_0}}=\prod_{\tau\in
J}\omega_{\tau}^{b_\tau}$, $\psi_2|_{I_{K_0}}=\prod_{\tau\notin
J}\omega_{\tau}^{b_\tau}$ with $2\leq b_\tau\leq p-2$ (note that for a fixed
$J$ it is \emph{not} always possible to do this, even after twisting -
indeed, up to twisting it is equivalent to $\rho$ being $J$-regular). In
this case we define canonical crystalline lifts $\psi_{1,J}$, $\psi_{2,J}$
of $\psi_1$, $\psi_2$, as in section \ref{2}. That is, we demand that
for some choice of a Frobenius element $\Frob_{K_0}\in G_{K_0}$, ${\psi}_{i,J}(\Frob_{K_0})$ is the Teichm\"{u}ller lift of $\psi_i(\Frob_{K_0})$, and that:
\begin{itemize}\item ${\psi}_{1,J}$ is crystalline, and the Hodge-Tate weight of ${\psi}_{1,J}$ with respect to $\tau$ is $b_\tau$ if $\tau\in J$, and $0$ if $\tau\notin J$.
\item ${\psi}_{2,J}$ is crystalline, and the Hodge-Tate weight of ${\psi}_{2,J}$ with respect to $\tau$ is $b_\tau$ if $\tau\notin J$, and $0$ if $\tau\in J$.\end{itemize}

The main result of
this section is
\begin{prop}\label{h1f}Under the above hypotheses, $\rho$ has a model of type $J$ if
and only if $\rho$ has a lift to a crystalline representation
$\bigl(\begin{smallmatrix}\psi_{1,J}&*\\0&\psi_{2,J}\end{smallmatrix}\bigr).$\end{prop}
\begin{proof}The idea of the proof is to express both the condition of
having a model of type $J$ and the condition of having a crystalline lift of
the prescribed type in terms of conditions on strongly divisible modules. In
fact, we already have a description of the general model of type $J$ in
terms of Breuil modules with descent data, and it is easy to write down the
general crystalline representation
$\bigl(\begin{smallmatrix}\psi_{1,J}&*\\0&\psi_{2,J}\end{smallmatrix}\bigr)$
in terms of Fontaine-Laffaille theory. The only difficulty comes in relating
the generic fibres of the Breuil modules to the generic fibres of the
Fontaine-Laffaille modules, as the image of the functors describing passage
to the generic fibre is in general too complicated to describe directly.
Fortunately, it is relatively easy to compare the two generic fibres we
obtain, without explicitly determining either.

Let $\M\in\BrMod^{k-1}_{dd}$ for some $k\in [2,p-1]$. Let $\hat{A}$ be
the filtered ring defined in section 2.1 of \cite{carusomodp}. There
is a contravariant functor $T_{st}^{*}$ from $\BrMod^{k-1}_{dd}$ to
the category of $E$-representations of $G_{K_{0}}$ given
by $$T_{st}^{*}(\M):=\Hom_{k[u]/u^{ep},\phi_{k-1},N,\Fil^\cdot}(\M,\hat{A})$$(where
compatibility with $\Fil^{\cdot}$ means that the image of $\M_{k-1}$
is contained in $\Fil^{k-1}\hat{A}$). The action of $G_{K_0}$ on
$T_{st}^{*}(\M)$ is given
by \[(gf)(x):=gf(\widehat{\overline{g}}^{-1}(x)),\]where
$\overline{g}$ is the image of $g$ in $\Gal(K/K_0)$, and the action of
$\Gal(K/K_0)$ on $\hat{A}=\hat{A}_K$ is defined in section 4.2 of
\cite{carusomodp}. For the compatibility of this definition with those
used in \cite{MR1681105}, \cite{MR2543474} and \cite{sav04}, see Lemma
3.3.1.2 of \cite{MR1621389}. This functor is exact and faithful, and
preserves dimension in the obvious sense. To see these properties, it
is enough to work with the category $\BrMod^{k-1}$ without descent
data, and it is also straightforward to see that it suffices to
consider the case $E=\Fp$. In this case, the fact that $T_{st}^*$ is
faithful is Corollary 2.3.3 of \cite{MR2543474}, and exactness follows
from Proposition 2.3.1 of \cite{carusomodp} and the duality explained
in section 2.1 of \cite{carusomodp}. The preservation of dimension is
Lemma 2.3.1.2 of \cite{MR1681105}.

We will see below that by Breuil's generalisation of
Fontaine-Laffaille theory (see \cite{MR1621389}) there are objects of
$\BrMod^{p-2}_{dd}$ which correspond via $T_{st}^{*}$ to the
reductions mod $\pi$ of crystalline representations with Hodge-Tate
weights in $[0,p-2]$. In order to compare the generic fibres of these
Breuil modules to those of finite flat group schemes with descent
data, we need to be able to compare elements of $\BrMod_{dd}^{1}$ and
$\BrMod_{dd}^{p-2}$. This is straightforward: it is easy to check that
there is a fully faithful functor from $\BrMod_{dd}^{1}$ to
$\BrMod_{dd}^{p-2}$, given by defining (for $\M\in\BrMod_{dd}^{1}$)
$\M_{p-2}:=u^{e(p-3)}\M_{1}$, $\phi_{p-2}(u^{e(p-3)}x)=\phi_{1}(x)$
for all $x\in\M_{1}$, and leaving the other structures unchanged. This
functor commutes with $T_{st}^{*}$.

Because we are now using the functor $T_{st}^{*}$ rather than
$T_{st,2}^{K_{0}}$, the form of the Breuil modules (and in particular
their descent data) corresponding to models of type $J$ under
$T_{st}^{*}$ is slightly different. We will simultaneously write it as
an element of $\BrMod_{dd}^1$ and $\BrMod_{dd}^{p-2}$ (making use of
the fully faithful functor of the previous paragraph), by specifying
$\mathcal{M}_1$, $\mathcal{M}_{p-2}$, $\phi_1$ and $\phi_{p-2}$. Explicitly, we see (recalling that
the operator $N$ is uniquely determined for an element of $\BrMod^1_{dd}$, so it suffices to check that
it satisfies $N(\mathcal{M})\subset u\mathcal{M}$ and the commutation
relations with $\phi_{1}$ and $\hat{g}$, which we will check below) from Lemma \ref{nicebasis} that $\rho$ has
a model of type $J$ if and only if there are $\lambda_{\tau_{i}}\in E$
with $\lambda_{\tau_{i}}=0$ if $\tau_{i+1}\notin J$, and elements $a$,
$b\in E^{\times}$ such that $\rho\cong T_{st}^{*}(\M)$, where

$$\M^{\tau_i}=E[u]/u^{ep}\cdot e_{\tau_i}+E[u]/u^{ep}\cdot f_{\tau_i}$$
$$\M_{1}^{\tau_i}=E[u]/u^{ep}\cdot u^{j_{\tau_i}}
e_{\tau_i}+E[u]/u^{ep}\cdot
(u^{e-j_{\tau_i}}f_{\tau_i}+\lambda_{\tau_i}u^{i_{\tau_i}}e_{\tau_i})$$
$$\M_{p-2}^{\tau_i}=E[u]/u^{ep}\cdot u^{(p-3)e+j_{\tau_i}} e_{\tau_i}+E[u]/u^{ep}\cdot (u^{(p-2)e-j_{\tau_i}}f_{\tau_i}+\lambda_{\tau_i}u^{(p-3)e+i_{\tau_i}}e_{\tau_i})$$
$$\phi_{1}(u^{j_{\tau_i}}
e_{\tau_i})=(a^{-1})_ie_{\tau_{i+1}} $$ $$\phi_{1}(u^{e-j_{\tau_i}}f_{\tau_i}+\lambda_{\tau_i}u^{i_{\tau_i}}e_{\tau_i})=(b^{-1})_if_{\tau_{i+1}}$$
$$\phi_{p-2}(u^{(p-3)e+j_{\tau_i}} e_{\tau_i})=(a^{-1})_ie_{\tau_{i+1}} $$ $$\phi_{p-2}(u^{(p-2)e-j_{\tau_i}}f_{\tau_i}+\lambda_{\tau_i}u^{(p-3)e+i_{\tau_i}}e_{\tau_i})=(b^{-1})_if_{\tau_{i+1}}$$
$$\hat{g}(e_{\tau_i})=\left(\left(\prod_{\sigma\notin
J}\omega_\sigma^{p-b_{\sigma}}\right)(g)\right)e_{\tau_i}$$
$$\hat{g}(f_{\tau_i})=\left(\left(\prod_{\sigma\in
J}\omega_\sigma^{p-b_{\sigma}}\right)(g)\right)f_{\tau_i}$$
$$N(e_{\tau_{i}})=0 $$
$$N(f_{\tau_{i}})=-\frac{(b)_{i-1}}{(a)_{i-1}}i_{\tau_{i-1}}\lambda_{\tau_{i-1}}u^{pi_{\tau_{i-1}}}e_{\tau_{i}} $$
where $\lambda_{\tau_i}\in
E$, with $\lambda_{\tau_i}=0$ if $\tau_{i+1}\notin J$, the $i_{\tau_i}$ are such
that $\M_{p-2}$ is Galois-stable and $0\leq i_{\tau_i}\leq e-1$, and
$$j_{\tau_i}=\left\{\begin{array}{c}
                 e\text{ if }\tau_{i+1}\in J \\
                 0\text{ if }\tau_{i+1}\notin J. \\
               \end{array}\right. $$
To see that $N(\mcal{M})\subset u\mcal{M}$, it is enough to check that
$i_{\tau_i}>0$ for all $i$. In fact, we claim that we have
$pi_{\tau_i}\ge e$ for all $i$. To see this, note that by definition
we have $i_{\tau_{i+1}}\equiv pi_{\tau_i}\pmod{e}$. If $pi_{\tau_i}<e$ for some $i$,
then this congruence forces $i_{\tau_{i+1}}=pi_{\tau_i}$. However, it is
easy to check that since $2\le b_{\tau_i}\le p-2$ for all $i$, no
$i_{\tau_i}$ is divisible by $p$  (for example,
by (6) below we have \[i_{\tau_i}\equiv
b_{\tau_i}(\delta_{J^c}(\tau_i)-\delta_J(\tau_i))-\delta_{J^c}(\tau_{i+1})\pmod{p},\]which
is never $0\pmod{p}$), so the claim follows. The compatibility of $N$
with $\hat{g}$ is evident from the definition of $i_{\tau_i}$.

To see that $u^eN(\mcal{M}_{1})\subset\mcal{M}_{1}$, and that
$\phi_{1}(u^eN(x))=N(\phi_{1}(x))$ for all $x\in\mcal{M}_{1}$,
we compute as follows (recalling that the Leibniz rule implies that
$N(u^ix)=u^iN(x)-iu^ix$): \[N(u^{j_{\tau_i}}e_{\tau_i})=-j_{\tau_i}u^{j_{\tau_i}}e_{\tau_i}\in\mcal{M}_{1},\]so
that \[\phi_{1}(u^eN(u^{j_{\tau_i}}e_{\tau_i}))=0=N((a^{-1})_ie_{\tau_{i+1}})=N(\phi_{1}(u^{j_{\tau_i}}e_{\tau_i})).\]Similarly,
we
have \begin{align*}N(u^{e-j_{\tau_i}}f_{\tau_i}+\lambda_{\tau_i}u^{i_{\tau_i}}e_{\tau_i})&=-\frac{(b)_{i-1}}{(a)_{i-1}}i_{\tau_{i-1}}\lambda_{\tau_{i-1}}u^{e-j_{\tau_i}+pi_{\tau_{i-1}}}e_{\tau_i}
  \\&- \left((e-j_{\tau_i})u^{e-j_{\tau_i}}f_{\tau_i}+i_{\tau_i}\lambda_{\tau_i}u^{i_{\tau_i}}e_{\tau_i}\right)  \end{align*} Recalling that if
$y\in\mcal{M}_{1}$ then $\phi_{1}(u^ey)=0$, we see that it is
enough to compute the right hand side modulo $\mcal{M}_{1}$. Since
$pi_{\tau_{i-1}}\ge e$, the exponent of $u$ in the first term on the
right hand side is at least $2e-j_{\tau_i}\ge j_{\tau_i}$,
so this term is contained in $\mcal{M}_{1}$. We thus see that modulo
$\mcal{M}_{1}$, the right hand side is congruent to
\[
  (e-j_{\tau_i}-i_{\tau_i})\lambda_{\tau_i}u^{i_{\tau_i}}e_{\tau_i}=-i_{\tau_i}\lambda_{\tau_i}u^{i_{\tau_i}}e_{\tau_i}
\](since if $\lambda_{\tau_i}\ne 0$, we have $\tau_{i+1}\in J$, so
that $j_{\tau_i}=e$). Then
\begin{align*}
  \phi_{1}(-u^ei_{\tau_i}\lambda_{\tau_i}u^{i_{\tau_i}}e_{\tau_i})&=\phi_{1}(-i_{\tau_i}\lambda_{\tau_i}u^{i_{\tau_i}}u^{j_{\tau_i}}e_{\tau_i})\\
  &=
  -i_{\tau_i}\lambda_{\tau_i}u^{pi_{\tau_i}}(a^{-1})_ie_{\tau_{i+1}}\\
  &=(b^{-1})_iN(f_{\tau_{i+1}})\\ &= N(\phi_{1}(u^{e-j_{\tau_i}}f_{\tau_i}+\lambda_{\tau_i}u^{i_{\tau_i}}e_{\tau_i}))
\end{align*}as required.

			It is an easy exercise to write down the
			reductions mod $p$ of the strongly divisible modules corresponding to
			crystalline representations
			$\bigl(\begin{smallmatrix}\psi_{1,J}&*\\0&\psi_{2,J}\end{smallmatrix}\bigr)$,
                        as we now explain. Firstly, we must recall one
                        of the main results of \cite{MR707328}. Let
                        $L$ be a finite extension of $\Qp$ with
                        residue field $E$. We
                        say that an admissible $\bigO_L$-lattice
                        is a finite free
                        $(\bigO_{K_0}\otimes_{\Z_p}\bigO_L)$-module $M$
                         together with a decreasing
                        filtration $\Fil^iM$ by
                        $\bigO_{K_0}$-direct summands and
                        $\phi$-linear, $\bigO_L$-linear maps
                        $\phi_i:\Fil^iM\to M$ for all $0\le i\le p-2$
                        such that
                        \begin{itemize}
                        \item $\Fil^0M=M$ and $\Fil^{p-1}M=0$.
                        \item For all $0\le i\le p-3$,
                          $\phi_i|_{\Fil^{i+1}M}=p\phi_{i+1}$.
                        \item $\sum_{i=0}^{p-2}\phi_i(\Fil^iM)=M$.
                        \end{itemize} There is an exact functor
                        $T^*_{cris}$ from the category of admissible
                        $\bigO_L$-lattices to the category of $G_{K_0}$-representations on free $\bigO_L$-lattices
                        defined
                        by \[T^*_{cris}(M)=\Hom_{\bigO_{K_0},\Fil^\cdot,\phi}(M,A_{cris}).\]This
                        gives an equivalence of categories between the
                        category of admissible $\bigO_L$-lattices and the
                        category of $G_{K_0}$-stable
                        $\bigO_L$-lattices in crystalline
                        $L$-representations in $G_{K_0}$ with all
                        Hodge-Tate weights in $[0,p-2]$.

In particular, one can easily write down the form of the rank one
$\bigO_L$-lattices corresponding to the characters $\psi_{1,J}$ and
$\psi_{2,J}$, and we must then compute the possible form of extensions
of these two lattices. As usual, we decompose $M$ as a direct sum of
$\bigO_L$-modules $M^{\tau_i}$. We obtain the following general
form:   \[M^{\tau_i}=\bigO_L E_{\tau_i}+\bigO_L
F_{\tau_i}\] \[\Fil^0M^{\tau_i}=M^{\tau_i}\] \[\Fil^{b_{\tau_i}+1}M^{\tau_i}=0\] \[\text{
  if }\tau_i\in J{\text, }\Fil^{j}M^{\tau_i}=\bigO_L F_{\tau_i}\text{ for
  all }1\le j\le b_{\tau_i}\]  \[\text{
  if }\tau_i\in J^c{\text, }\Fil^{j}M^{\tau_i}=\bigO_L E_{\tau_i}\text{ for
  all }1\le j\le b_{\tau_i}\]  \[\text{
  if }\tau_i\in J{\text, }\phi_0(E_{\tau_i})=(\tilde{a}^{-1})_iE_{\tau_{i+1}}\text{ and }\phi_{b_{\tau_i}}(F_{\tau_i})=(\tilde{b}^{-1})_i(F_{\tau_{i+1}}-\lambda'_{\tau_i}E_{\tau_{i+1}}) \]  \[\text{
  if }\tau_i\notin J{\text,
}\phi_{b_{\tau_i}}(E_{\tau_i})=(\tilde{a}^{-1})_iE_{\tau_{i+1}}\text{
  and
}\phi_{0}(F_{\tau_i})=(\tilde{b}^{-1})_i(F_{\tau_{i+1}}-\lambda'_{\tau_i}E_{\tau_{i+1}}) \]
where $\tilde{a}$, $\tilde{b}\in\bigO_L^\times$, $\lambda'_{\tau_i}\in\bigO_L$, and $\lambda'_{\tau_i}=0$ if
$\tau_{i+1}\notin J$. To see this, note that the form of the
filtration is easily deduced from the relationship between the
filtration of a Fontaine-Laffaille module, and the Hodge-Tate weights
of the corresponding Galois representation, and the form of the
Frobenius action on the $E_{\tau_i}$ is also determined. To see that
we can arrange the Frobenius action as claimed, suppose firstly that
$J\neq\emptyset$, and choose $\tau_i\in J$. Then $F_{\tau_i}$ is
determined (by the form of the filtration)
up to an element of $\bigO_L^\times$, and we fix a choice of $F_{\tau_i}$. If
$\tau_{i+1}\notin J$, we can simply define
$F_{\tau_{i+1}}=(\tilde{b})_i\phi_{b_{\tau_i}}(F_{\tau_i})$. If
$\tau_{i+1}\in J$, then there is a unique
$\lambda_{\tau_i}'\in\bigO_L$ such
that \[\Fil^{b_{\tau_{i+1}}}M^{\tau_{i+1}}=\bigO_L((\tilde{b})_i\phi_{b_{\tau_i}}(F_{\tau_i})+\lambda'_{\tau_i}E_{\tau_{i+1}}),\]
and we
set \[F_{\tau_{i+1}}=(\tilde{b})_i\phi_{b_{\tau_i}}(F_{\tau_i})+\lambda'_{\tau_i}E_{\tau_{i+1}}.\]We
can then continue in the same fashion, defining $F_{\tau_{i+2}}$ and so
on, and
the fact that $\tau_i\in J$ gives us the freedom to choose
$\lambda'_{\tau_{i-1}}$ so
that \[\phi_{\delta_J(\tau_{i-1})b_{\tau_{i-1}}}(F_{\tau_{i-1}})=(\tilde{b}^{-1})_{i-1}(F_{\tau_{i}}-\lambda'_{\tau_{i-1}}E_{\tau_{i}})
.\]The case $J=\emptyset$ is similar, except that one may need to
modify the initial choice of $F_{\tau_i}$; the argument is very
similar to that used in the case $J=S$ in the proof of Lemma
\ref{nicebasis}. In this case one also needs to use the fact that
$\tilde{b}^{-1}-\tilde{a}^{-1}p^{\sum_{\tau\in
    S}b_\tau}\in\bigO_L^\times$, which holds as $\sum_{\tau\in S}b_\tau>0$.

Breuil's generalisation of Fontaine-Laffaille theory (\cite{MR1621389}) allows us to
reduce these modules mod $\pi_L$ and obtain the corresponding elements
of the category $K_0-\BrMod_{dd,K_0}^{p-2}$ of Breuil modules with
descent data for the case $K=K_0$ (so that the descent data is trivial,
as the group $\Gal(K/K_0)$ is trivial). We find that they are of the
form: 	$$\mathcal{Q}^{\tau_i}=E[u]/u^{p}\cdot E_{\tau_i}+E[u]/u^{p}\cdot F_{\tau_i}$$
			$$\mathcal{Q}_{p-2}^{\tau_i}=E[u]/u^{p}\cdot u^{(p-2-b_{\tau_{i}}\delta_{J^{c}}(\tau_{i}))} E_{\tau_i}+E[u]/u^{p}\cdot u^{(p-2-b_{\tau_{i}}\delta_{J}(\tau_{i}))}F_{\tau_i}$$
			$$\phi_{p-2}(u^{(p-2-b_{\tau_{i}}\delta_{J^{c}}(\tau_{i}))} E_{\tau_i})=(a^{-1})_iE_{\tau_{i+1}} $$ $$\phi_{p-2}(u^{(p-2-b_{\tau_{i}}\delta_{J}(\tau_{i}))}F_{\tau_i})=(b^{-1})_i(F_{\tau_{i+1}}-\lambda_{\tau_{i}}'E_{\tau_{i+1}})$$
					$$N(E_{\tau_{i}})=0 $$
			$$N(F_{\tau_{i}})=0 $$
			where $\lambda_{\tau_i}'\in
			E$, with $\lambda'_{\tau_i}=0$ if $\tau_{i+1}\notin J$.

                        Of course, we wish to know the corresponding
                        objects of $\BrMod_{dd}^{p-2}$. This is
                        straightforward: by Proposition 4.2.2 of
                        \cite{carusomodp}, and the discussion
                        preceding and following it, we see that we can
                        obtain the requisite modules by simply taking
                        the extension of scalars $k[u]/u^{p}\to
                        k[u]/u^{ep}$ given by $u\mapsto u^e$, and
                        allowing $\Gal(K/K_0)$ to act via its action
                        on $k[u]/u^{ep}$.  We obtain the
                        following general form:
			
			$$\mathcal{N}^{\tau_i}=E[u]/u^{ep}\cdot E_{\tau_i}+E[u]/u^{ep}\cdot F_{\tau_i}$$
			$$\mathcal{N}_{p-2}^{\tau_i}=E[u]/u^{ep}\cdot u^{e(p-2-b_{\tau_{i}}\delta_{J^{c}}(\tau_{i}))} E_{\tau_i}+E[u]/u^{ep}\cdot u^{e(p-2-b_{\tau_{i}}\delta_{J}(\tau_{i}))}F_{\tau_i}$$
			$$\phi_{p-2}(u^{e(p-2-b_{\tau_{i}}\delta_{J^{c}}(\tau_{i}))} E_{\tau_i})=(a^{-1})_iE_{\tau_{i+1}} $$ $$\phi_{p-2}(u^{e(p-2-b_{\tau_{i}}\delta_{J}(\tau_{i}))}F_{\tau_i})=(b^{-1})_i(F_{\tau_{i+1}}-\lambda_{\tau_{i}}'E_{\tau_{i+1}})$$
			$$\hat{g}(E_{\tau_i})=E_{\tau_i}$$
			$$\hat{g}(F_{\tau_i})=F_{\tau_i}$$
			$$N(E_{\tau_{i}})=0 $$
			$$N(F_{\tau_{i}})=0 $$
			where $\lambda_{\tau_i}'\in
			E$, with $\lambda'_{\tau_i}=0$ if $\tau_{i+1}\notin J$.
We claim that if for each $i$ we have			\begin{equation}\label{eqn:341}
\lambda_{\tau_{i}}(b)_{i}=\lambda'_{\tau_{i}}(a)_{i}
			\end{equation} then $T_{st}^{*}(\M)\cong T_{st}^{*}(\mathcal{N})$. This is of course enough to demonstrate the proposition, as given any set of $\lambda_{\tau_{i}}$ (respectively $\lambda_{\tau_{i}}'$) such that $\lambda_{\tau_i}=0$ (respectively $\lambda'_{\tau_i}=0$) if $\tau_{i+1}\notin J$, we may choose a set of $\lambda_{\tau_{i}}'$ (respectively $\lambda_{\tau_{i}}$) so that (\ref{eqn:341}) holds.

Assume now that (\ref{eqn:341}) holds. Note that we may write both $\M$ and $\mathcal{N}$ as extensions 
$$0\to\M''\to\M\to\M'\to 0 $$
$$0\to\mathcal{N}''\to\mathcal{N}\to\mathcal{N}'\to 0 $$with $T_{st}^{*}(\M'')\cong T_{st}^{*}(\mathcal{N}'')\cong\psi_{2}$, $T_{st}^{*}(\M')\cong T_{st}^{*}(\mathcal{N}')\cong\psi_{1}$.

 To prove that $T_{st}^{*}(\M)\cong T_{st}^{*}(\mathcal{N})$, we will construct a commutative diagram 
$$\xymatrix{0\ar[r] & \M''\ar[r]\ar[d]^{f_{\mathcal{M}''}}
  &\M\ar[r]\ar[d]^{f_{\mathcal{M}}} &
  \M'\ar[r]\ar[d]^{f_{\mathcal{M}'}} & 0 \\ 0\ar[r] &
  \mathcal{P}''\ar[r] &\mathcal{P}\ar[r] & \mathcal{P}'\ar[r] & 0\\
  0\ar[r] & \mathcal{N}''\ar[u]_{f_{\mathcal{N}''}}\ar[r]
  &\mathcal{N}\ar[r]\ar[u]_{f_{\mathcal{N}}} &
  \mathcal{N}'\ar[r]\ar[u]_{f_{\mathcal{N}'}} & 0}$$such that each of
$T_{st}^{*}(f_{\mathcal{M}''})$, $T_{st}^{*}(f_{\mathcal{M}'})$,
$T_{st}^{*}(f_{\mathcal{N}''})$ and $T_{st}^{*}(f_{\mathcal{N}'})$ are
isomorphisms. From the five lemma it then follows that
$T_{st}^{*}(f_{\mathcal{M}})$ and $T_{st}^{*}(f_{\mathcal{N}})$ are
isomorphisms, and we will be done.

In fact, we take

$$\mathcal{P}^{\tau_i}=E[u]/u^{ep}\cdot e'_{\tau_i}+E[u]/u^{ep}\cdot f'_{\tau_i}$$
$$\mathcal{P}_{p-2}^{\tau_i}=E[u]/u^{ep}\cdot u^{n_{\tau_{i}}} e'_{\tau_i}+E[u]/u^{ep}\cdot (u^{n'_{\tau_{i}}}f'_{\tau_i}+\lambda_{\tau_i}u^{n_{\tau_{i}}-\beta_{\tau_{i+1}}}e'_{\tau_i})$$
$$\phi_{p-2}(u^{n_{\tau_{i}}} e'_{\tau_i} )=(a^{-1})_ie'_{\tau_{i+1}} $$ $$\phi_{p-2}(u^{n'_{\tau_{i}}}f'_{\tau_i}+\lambda_{\tau_i}u^{n_{\tau_{i}}-\beta_{\tau_{i+1}}}e'_{\tau_i})=(b^{-1})_if'_{\tau_{i+1}}$$
$$\hat{g}(e'_{\tau_i})=\nu_{1,\tau_{i}}(g)e'_{\tau_i}$$
$$\hat{g}(f'_{\tau_i})=\nu_{2,\tau_{i}}(g)f'_{\tau_i}$$
$$N(e'_{\tau_{i}})=0 $$
$$N(f'_{\tau_{i}})=-\frac{(b)_{i-1}}{(a)_{i-1}}i_{\tau_{i-1}}\lambda_{\tau_{i-1}}u^{pi_{\tau_{i-1}}-p \alpha_{\tau_{i}}}e'_{\tau_{i}} $$where $$\alpha_{\tau_{i}}=\sum_{j=0}^{r-1}p^{r-1-j}\left(b_{\tau_{i+j}}\delta_{J^{c}}(\tau_{i+j})-\delta_{J^{c}}(\tau_{i+j+1})\right) $$
$$\beta_{\tau_{i}}=\sum_{j=0}^{r-1}p^{r-1-j}\left(b_{\tau_{i+j}}\delta_{J}(\tau_{i+j})-\delta_{J}(\tau_{i+j+1})\right) $$ $$\nu_{1,\tau_{i}}(g)=\left\{\begin{array}{ll}\prod_{\sigma\notin
J}\omega_\sigma^{p-b_{\sigma}}(g)&\text{ if }\tau_i\notin J \\ 1 &\text{ if }\tau_i\in J \end{array}\right.$$ $$\nu_{2,\tau_{i}}(g)=\left\{\begin{array}{ll}\prod_{\sigma\in
J}\omega_\sigma^{p-b_{\sigma}}(g)&\text{ if }\tau_i\in J \\ 1 &\text{ if }\tau_i\notin J \end{array}\right.$$ $$n_{\tau_{i}}= (p-2-\delta_{J^{c}}(\tau_{i})b_{\tau_{i}})e+p\delta_{J^{c}}(\tau_{i})\alpha_{\tau_{i}}-\delta_{J^{c}}(\tau_{i+1})\alpha_{\tau_{i+1}}$$
$$n'_{\tau_{i}}= (p-2-\delta_{J}(\tau_{i})b_{\tau_{i}})e+p\delta_{J}(\tau_{i})\beta_{\tau_{i}}-\delta_{J}(\tau_{i+1})\beta_{\tau_{i+1}}$$We then define $f_{\M}$ and $f_{\mathcal{N}}$ by $$f_{\M}(e_{\tau_{i}})=u^{-p\alpha_{\tau_{i}}\delta_{J}(\tau_{i})}e_{\tau_{i}}' $$
$$f_{\M}(f_{\tau_{i}})=u^{-p\beta_{\tau_{i}}\delta_{J^{c}}(\tau_{i})}f_{\tau_{i}}' $$
$$f_{\mathcal{N}}(E_{\tau_{i}})=u^{p\alpha_{\tau_{i}}\delta_{J^{c}}(\tau_{i})}e_{\tau_{i}}' $$
$$f_{\mathcal{N}}(F_{\tau_{i}})=u^{p\beta_{\tau_{i}}\delta_{J}(\tau_{i})}f_{\tau_{i}}' $$We define $\mathcal{P}'$ to be the submodule generated by the $e'_{\tau_{i}}$, and $\mathcal{P}''$ to be the quotient obtained by $e'_{\tau_{i}}\mapsto 0$. The remaining maps are then defined by the commutativity of the diagram.

Before we verify that this construction behaves as claimed, we pause to
record a number of useful identities and inequalities.
\begin{enumerate}
\item If $\tau_{i+1}\notin J$, then
  $\lambda_{\tau_i}=\lambda'_{\tau_i}=0$ by definition.
\item  $$p\alpha_{\tau_{i}}-\alpha_{\tau_{i+1}}=e(b_{\tau_{i}}\delta_{J^{c}}(\tau_{i})-\delta_{J^{c}}(\tau_{i+1})), $$
$$p\beta_{\tau_{i}}-\beta_{\tau_{i+1}}=e(b_{\tau_{i}}\delta_{J}(\tau_{i})-\delta_{J}(\tau_{i+1})). $$These
both follow immediately from the definitions of $\alpha_{\tau_i}$, $\beta_{\tau_i}$.
\item \[n_{\tau_i}=\alpha_{\tau_{i+1}}\delta_J(\tau_{i+1})-p\alpha_{\tau_i}\delta_J(\tau_i)+e(p-3)+e\delta_J(\tau_{i+1}),\] \[n'_{\tau_i}=\beta_{\tau_{i+1}}\delta_{J^c}(\tau_{i+1})-p\beta_{\tau_i}\delta_{J^c}(\tau_i)+e(p-3)+e\delta_{J^c}(\tau_{i+1}).\]These
  both follow from the definitions of $n_{\tau_i}$, $n_{\tau_i}'$ and property (2) above.
\item We have $\tau_{i}\in J$ if and only if $\beta_{\tau_{i}}>0$ if
  and only if $\alpha_{\tau_{i}}\leq 0$. To see this, note that from
  the definition, the sign of $\alpha_{\tau_{i}}$ is determined by the
  sign of the first non-zero term in the sum (this uses that $2\leq
  b_{\tau_{j}}\leq p-2$). If $\tau_{i}\notin J$ then the first term is
  positive, and thus so is the whole sum. If $\tau_{i}\in J$ then
  either every term in the sum is zero, or the first non-zero term
  must be negative. A similar analysis applies to the sign of
  $\beta_{\tau_{i}}$.
\item  $$-e/(p-1)<\alpha_{\tau_{i}},\ \beta_{\tau_{i}}<
 e(p-2)/(p-1).$$ This is immediate from the definitions, and the fact that $2\leq b_{\tau_{j}}\leq p-2$ for all
 $j$.
\item \[i_{\tau_{i-1}}=\alpha_{\tau_i}-\beta_{\tau_i}+e\delta_J(\tau_i).\] It
  follows straightforwardly from the forms of the $\hat{g}$-actions
  that the two side are congruent modulo $e$, so it suffices to check
  that the right hand side is an element of $[0,e-1]$. This follows
  from points (4) and (5).
\item \[(p-2)e\ge n_{\tau_i},\ n_{\tau_i}\ge 0.\] We demonstrate these
  inequalities for $n_{\tau_i}$, the argument for $n'_{\tau_i}$ being
  formally identical after exchanging $\alpha_{\tau_j}$ and
  $\beta_{\tau_j}$, $J$ and $J^c$. We examine 4 cases in turn. If
  $\tau_{i}\in J$ and $\tau_{i+1}\in J$, then $n_{\tau_i}=(p-2)e$ and there
  is nothing to prove. If $\tau_i\in J$ and $\tau_{i+1}\notin J$, then
  $n_{\tau_i}=(p-2)e-\alpha_{\tau_{i+1}}$, and the inequalities follow
  from points (4) and (5) above. If $\tau_i\notin J$ and
  $\tau_{i+1}\in J$ then by point (3) above we have
  $n_{\tau_i}=(p-2)e+\alpha_{\tau_{i+1}}$, and the
  inequalities follow from (4) and (5). Finally, if $\tau_i\notin J$
  and $\tau_{i+1}\notin J$, then
  $n_{\tau_i}=(p-2-b_{\tau_i})e+p\alpha_{\tau_i}-\alpha_{\tau_{i+1}}=e(p-3)$
  by (2).
\item If $\tau_{i+1}\in J$, we
  have \[n_{\tau_i}-\beta_{\tau_{i+1}}=e(p-3)-p\alpha_{\tau_i}\delta_J(\tau_i)+i_{\tau_i}.\]This
  follows from (3) and (6) above.
\item If $\tau_{i+1}\in J$, then
  \[n_{\tau_i}-\beta_{\tau_{i+1}}\equiv n_{\tau_i}'+
  i_{\tau_i}\pmod{p}.\]This follows from (3) and (8).
\item If $\tau_{i+1}\in J$, then $n_{\tau_i}\ge
  \beta_{\tau_{i+1}}$. This follows from (4) and (8).
\item If $\tau_{i+1}\in J$, then \[n'_{\tau_i}+\beta_{\tau_{i+1}}\le
  e(p-2).\]From (2) and (3), we
  obtain \[n'_{\tau_{i}}+\beta_{\tau_{i+1}}=e(p-2)+\delta_J(\tau_i)(p\beta_{\tau_i}-eb_{\tau_i}),\]so
  we must check that if $\tau_i\in J$, then
  $p\beta_{\tau_i}-eb_{\tau_i}\le 0$. But by the definition of
  $\beta_{\tau_i}$, if $\tau_i\in J$ then we
  have \[p\beta_{\tau_i}-eb_{\tau_i}=-p^r+b_{\tau_i}+\sum_{j=1}^{r-1}p^{r-j}\left(b_{\tau_{i+j}}\delta_{J}(\tau_{i+j})-\delta_{J}(\tau_{i+j+1})\right)\]and
  the result follows as $2\le b_{\tau_j}\le p-2$ for all $j$.
\item If $\tau_i\in J$, then
  \[n_{\tau_i}'+pi_{\tau_{i-1}}-p\alpha_{\tau_i}\ge n_{\tau_i}.\]To
  see this, by (2), (3) and (6) we
  have that if $\tau_{i+1}\in J$,
  then  \begin{align*}n_{\tau_i}'+pi_{\tau_{i-1}}-p\alpha_{\tau_i}-n_{\tau_i}&=(p-1)e-p\beta_{\tau_i}\\
    &\ge\left((p-1)-\frac{p(p-2)}{p-1}\right)e\\&\ge 0\end{align*}by
  (5). If on the other hand $\tau_{i+1}\notin J$, we find that  \begin{align*}n_{\tau_i}'+pi_{\tau_{i-1}}-p\alpha_{\tau_i}-n_{\tau_i}&=(p+1-b_{\tau_i})e+p\alpha_{\tau_i}\\
   &\ge 3e+p\alpha_{\tau_i} \\&\ge\left(3-\frac{p}{p-1}\right)e\\&\ge 0\end{align*} by (5). \end{enumerate}

We now verify that $\mcal{P}$ is indeed an object of
$\BrMod_{dd}^{p-2}$. 
\begin{itemize}
\item To see that we have defined a $(k\otimes_{\Fp}
  E)[u])/u^{ep}$-module, we must check that all of the exponents of
  $u$ in the definition are nonnegative; so, we need to check the
  inequalities $n_{\tau_i}\ge 0$, $n'_{\tau_i}\ge 0$, and if
  $\lambda_{\tau_i}\ne 0$ we need to verify that
  $n_{\tau_i}\ge\beta_{\tau_{i+1}}$. These follow from (1), (7) and
  (10) above.
\item To see that $u^{e(p-2)}\mcal{P}\subset\mcal{P}_{p-2}$, we need to
  verify that $(p-2)e\ge n_{\tau_i}$, that $(p-2)e\ge n'_{\tau_i}$,
  and if $\lambda_{\tau_i}\ne 0$ then $(p-2)e\ge
  n'_{\tau_i}+\beta_{\tau_{i+1}}$. These follow from (1), (7) and (11).
\item To see that $N(\mcal{P})\subset u\mcal{P}$, we need to check
  that if $\lambda_{\tau_{i-1}}\ne 0$ then
  $i_{\tau_{i-1}}>\alpha_{\tau_i}$. This follows from (1), (5) and (6).
\item To see that $u^e N(\mcal{P}_{p-2})\subset \mcal{P}_{p-2}$, we
  note by the Leibniz rule we
  have \[N(u^{n_{\tau_i}}e'_{\tau_i})=-n_{\tau_i}u^{n_{\tau_i}}e'_{\tau_i}\in
  \mcal{P}_{p-2}\]and \begin{align*}N(u^{n'_{\tau_{i}}}f'_{\tau_i}+\lambda_{\tau_i}u^{n_{\tau_{i}}-\beta_{\tau_{i+1}}}e'_{\tau_i})&=
 -\frac{(b)_{i-1}}{(a)_{i-1}}i_{\tau_{i-1}}\lambda_{\tau_{i-1}} u^{n'_{\tau_{i}}+pi_{\tau_{i-1}}-p
      \alpha_{\tau_{i}}}e'_{\tau_{i}}
  \\&-n'_{\tau_i}u^{n'_{\tau_i}}f'_{\tau_i}-\lambda_{\tau_i}(n_{\tau_{i}}-\beta_{\tau_{i+1}})u^{n_{\tau_{i}}-\beta_{\tau_{i+1}}}e'_{\tau_i}.\end{align*}By
(1) and
(12), the first time on the right hand side is contained in
$\mcal{P}_{p-2}$. By (9), the remaining terms are equal
to \[-n'_{\tau_i}(u^{n'_{\tau_i}}f'_{\tau_i}+\lambda_{\tau_i}u^{n_{\tau_{i}}-\beta_{\tau_{i+1}}}e'_{\tau_i})-\lambda_{\tau_i}i_{\tau_{i}}u^{n_{\tau_{i}}-\beta_{\tau_{i+1}}}e'_{\tau_i}.\]The
first term is contained in $\mcal{P}_{p-2}$ by definition, and if we
multiply the second term by $u^e$, we
obtain \[-\lambda_{\tau_i}i_{\tau_{i}}u^{n_{\tau_{i}}+e-\beta_{\tau_{i+1}}}e'_{\tau_i},\]which
is contained in $\mcal{P}_{p-2}$ by (5).
\item To see that $\phi_{p-2}(u^eN(x))=N(\phi_{p-2}(x))$ for all $x\in
  \mcal{P}_{p-2}$, we recall that $\phi_{p-2}(u^ey)=0$ if
  $y\in\mcal{P}_{p-2}$. Thus \[\phi_{p-2}(u^eN(u^{n_{\tau_i}}e'_{\tau_i}))=0=N((a^{-1})_ie'_{\tau_{i+1}})=N(\phi_{p-2}(u^{n_{\tau_i}}e'_{\tau_i})).\]
  We also have, using (1), (6) and the calculation of the previous
  bullet point,
  \begin{align*}\phi_{p-2}(u^eN(u^{n'_{\tau_{i}}}f'_{\tau_i}+\lambda_{\tau_i}u^{n_{\tau_{i}}-\beta_{\tau_{i+1}}}e'_{\tau_i}))&=\phi_{p-2}(-\lambda_{\tau_i}i_{\tau_{i}}u^{n_{\tau_{i}}+e-\beta_{\tau_{i+1}}}e'_{\tau_i})\\
    &=-\lambda_{\tau_i}i_{\tau_i}u^{p(e-\beta_{i+1})}(a^{-1})_{i}e'_{\tau_{i+1}}\\
    &=-\lambda_{\tau_i}i_{\tau_i}u^{p(i_{\tau_i}-\alpha_{\tau_{i+1}})}(a^{-1})_{i}e'_{\tau_{i+1}}\\ &=(b^{-1})_{i}N(f'_{\tau_{i+1}})\\&=N(\phi_{p-2}(u^{n'_{\tau_{i}}}f'_{\tau_i}+\lambda_{\tau_i}u^{n_{\tau_{i}}-\beta_{\tau_{i+1}}}e'_{\tau_i})).\end{align*}

\item That $\mcal{P}_{p-2}$ is $\hat{g}$-stable follows directly from
  the definitions of  $\beta_{\tau_i}$, $n_{\tau_i}$ and $n'_{\tau_i}$.
\item That the action of $\hat{g}$ commutes with
  $\phi_{p-2}$ follows from the definition of $n_{\tau_i}$ and $n'_{\tau_i}$.
\item That the action of $\hat{g}$ commutes with
  $N$ follows from (6).
\end{itemize}

We now verify the claimed properties of $f_{\M}$ and $f_{\mcal{N}}$.
\begin{itemize}
\item In order that the maps $f_{\M}$ and $f_{\mathcal{N}}$ be
  defined, it is necessary that the exponents of $u$ in their
  definition be non-negative. This follows from (4).
\item To see that $f_{\M}(\M_{p-2})\subset\mcal{P}_{p-2}$ and
  $f_{\mcal{N}}(\mcal{N}_{p-2})\subset\mcal{P}_{p-2}$, we compute as
  follows. \begin{align*}f_{\mcal{M}}(u^{(p-3)e+j_{\tau_i}}e_{\tau_i})&=u^{(p-3)e+j_{\tau_i}-n_{\tau_i}-p\alpha_{\tau_i}\delta_J(\tau_i)}u^{n_{\tau_i}}e_{\tau_i}'\\&=u^{-\delta_J(\tau_{i+1})\alpha_{\tau_{i+1}}}(u^{n_{\tau_i}}e_{\tau_i}')\end{align*}by
  (3) and the definition of $j_{\tau_i}$. Similarly, by using (1), (3) and (6), we find
  that \[f_{\mcal{M}}(u^{(p-2)e-j_{\tau_i}}f_{\tau_i}+\lambda_{\tau_i}u^{(p-3)e+i_{\tau_i}}e_{\tau_i})=u^{-\beta_{\tau_{i+1}}\delta_{J^c}(\tau_{i+1})}(u^{n'_{\tau_i}}f_{\tau_i}'+\lambda_{\tau_i}u^{n_{\tau_i}-\beta_{\tau_{i+1}}}e_{\tau_i}').\]In
  the same
  way, using (1) and the definitions of $n_{\tau_i}$ and $n'_{\tau_i}$, \[f_{\mcal{N}}(u^{e(p-2-b_{\tau_i}\delta_{J^c}(\tau_i))}E_{\tau_i})=u^{\delta_{J^c}(\tau_{i+1})\alpha_{\tau_{i+1}}}(u^{n_{\tau_i}}e_{\tau_i}'),\] \[f_{\mcal{N}}(u^{e(p-2-b_{\tau_i}\delta_J(\tau_i))}F_{\tau_i})=u^{\delta_{J}(\tau_{i+1})\beta_{\tau_{i+1}}}(u^{n'_{\tau_i}}f'_{\tau_i}+\lambda_{\tau_i}u^{n_{\tau_i}-\beta_{\tau_{i+1}}}e'_{\tau_{i}})-\lambda_{\tau_i}u^{n_{\tau_i}}e'_{\tau_i}.  \]The
  result then follows from (4).

\item To check that $f_{\mcal{M}}$ and $f_{\mcal{N}}$ commute with
  $\phi_{p-2}$, we again compute directly. We
  have \begin{align*}f_{\mcal{M}}(\phi_{p-2}(u^{(p-3)e+j_{\tau_i}}e_{\tau_i}))&=f_{\mcal{M}}((a^{-1})_ie_{\tau_{i+1}})\\&=(a^{-1})_iu^{-p\alpha_{\tau_{i+1}}\delta_J(\tau_{i+1})}e'_{\tau_{i+1}},\end{align*}while
  \begin{align*}\phi_{p-2}(f_{\mcal{M}}(u^{(p-3)e+j_{\tau_i}}e_{\tau_i}))&=
  \phi_{p-2}(u^{-\delta_J(\tau_{i+1})\alpha_{\tau_{i+1}}}(u^{n_{\tau_i}}e_{\tau_i}'))\\&=(a^{-1})_iu^{-p\alpha_{\tau_{i+1}}\delta_J(\tau_{i+1})}e'_{\tau_{i+1}}.\end{align*}Similarly,
  we
  find \begin{align*}f_{\mcal{M}}(\phi_{p-2}(u^{(p-2)e-j_{\tau_i}}f_{\tau_i}+\lambda_{\tau_i}u^{(p-3)e+i_{\tau_i}}e_{\tau_i}))&=f_{\mcal{M}}((b^{-1})_if_{\tau_{i+1}})\\&=(b^{-1})_iu^{-p\beta_{\tau_{i+1}}\delta_{J^c}(\tau_{i+1})}f'_{\tau_{i+1}},\end{align*} \begin{align*}\phi_{p-2}(f_{\mcal{M}}(u^{(p-2)e-j_{\tau_i}}f_{\tau_i}+\lambda_{\tau_i}u^{(p-3)e+i_{\tau_i}}e_{\tau_i})))&=
  \phi_{p-2}(u^{-\beta_{\tau_{i+1}}\delta_{J^c}(\tau_{i+1})}(u^{n'_{\tau_i}}f_{\tau_i}'+\lambda_{\tau_i}u^{n_{\tau_i}-\beta_{\tau_{i+1}}}e_{\tau_i}')))\\&=(b^{-1})_iu^{-p\beta_{\tau_{i+1}}\delta_{J^c}(\tau_{i+1})}f'_{\tau_{i+1}},\end{align*} \begin{align*}f_{\mcal{N}}(\phi_{p-2}(u^{e(p-2-b_{\tau_i}\delta_{J^c}(\tau_i))}E_{\tau_i}))&=f_{\mcal{N}}((a^{-1})_iE_{\tau_{i+1}})\\&=(a^{-1})_iu^{p\delta_{J^c}(\tau_{i+1})\alpha_{\tau_{i+1}}}e'_{\tau_{i+1}}, \end{align*} \begin{align*}\phi_{p-2}(f_{\mcal{N}}(u^{e(p-2-b_{\tau_i}\delta_{J^c}(\tau_i))}E_{\tau_i}))&=
    \phi_{p-2}(u^{\delta_{J^c}(\tau_{i+1})\alpha_{\tau_{i+1}}}(u^{n_{\tau_i}}e_{\tau_i}')\\&=(a^{-1})_iu^{p\delta_{J^c}(\tau_{i+1})\alpha_{\tau_{i+1}}}e'_{\tau_{i+1}},\end{align*} \begin{align*}f_{\mcal{N}}(\phi_{p-2}(u^{e(p-2-b_{\tau_i}\delta_J(\tau_i))}F_{\tau_i}))&=f_{\mcal{N}}((b^{-1})_i(F_{\tau_{i+1}}-\lambda'_{\tau_{i}}E_{\tau_{i+1}}))\\&=(b^{-1})_i(u^{p\delta_{J}(\tau_{i+1})\beta_{\tau_{i+1}}}f'_{\tau_{i+1}}-\lambda'_{\tau_i}u^{p\delta_{J^c}(\tau_{i+1})\alpha_{\tau_{i+1}}}e'_{\tau_{i+1}}), \end{align*} \begin{align*}
      \phi_{p-2}(f_{\mcal{N}}(u^{e(p-2-b_{\tau_i}\delta_J(\tau_i))}F_{\tau_i}))&=\phi_{p-2}(u^{\delta_{J}(\tau_{i+1})\beta_{\tau_{i+1}}}(u^{n'_{\tau_i}}f'_{\tau_i}+\lambda_{\tau_i}u^{n_{\tau_i}-\beta_{\tau_{i+1}}}e'_{\tau_{i}})-\lambda_{\tau_i}u^{n_{\tau_i}}e'_{\tau_i})\\&=u^{p\delta_{J}(\tau_{i+1})\beta_{\tau_{i+1}}}(b^{-1})_if'_{\tau_{i+1}}-\lambda_{\tau_i}(a^{-1})_ie'_{\tau_{i+1}}.\end{align*}
      The result follows, because
      $\lambda_{\tau_i}(a^{-1})_i=\lambda'_{\tau_i}(b^{-1})_i$ by (~(\ref{eqn:341})), and if
      $\lambda_{\tau_i}\ne 0$ then $\delta_{J^c}(\tau_{i+1})=0$ by (1).
\item To check that $f_{\mcal{M}}$ and $f_{\mcal{N}}$ commute with
  $N$, we again compute directly. We
  have \begin{align*}N(f_{\mcal{M}}(e_{\tau_i})&=N(u^{-p\alpha_{\tau_i}\delta_{J}(\tau_i)}e'_{\tau_i})\\&=-p\alpha_{\tau_i}\delta_{J}(\tau_i)u^{-p\alpha_{\tau_i}\delta_{J}(\tau_i)}e'_{\tau_i}\\&=0\\&=f_{\mcal{M}}(N(e_{\tau_i})). \end{align*} Similarly,\begin{align*}
  N(f_{\mcal{M}}(f_{\tau_i}))&=N(u^{-p\beta_{\tau_i}\delta_{J^c}(\tau_i)}f'_{\tau_i})\\&=u^{-p\beta_{\tau_i}\delta_{J^c}(\tau_i)}N(f'_{\tau_i})\\&=
  -\frac{(b)_{i-1}}{(a)_{i-1}}i_{\tau_{i-1}}\lambda_{\tau_{i-1}}u^{pi_{\tau_{i-1}}-p
    \alpha_{\tau_{i}}-p\beta_{\tau_i}\delta_{J^c}(\tau_i)}e'_{\tau_{i}},\end{align*}while \begin{align*}f_{\mcal{M}}(N(f_{\tau_i}))&=f_{\mcal{M}}\left(-\frac{(b)_{i-1}}{(a)_{i-1}}i_{\tau_{i-1}}\lambda_{\tau_{i-1}}u^{pi_{\tau_{i-1}}}e_{\tau_{i}}
  \right)\\&=-\frac{(b)_{i-1}}{(a)_{i-1}}i_{\tau_{i-1}}\lambda_{\tau_{i-1}}u^{pi_{\tau_{i-1}}-p\alpha_{\tau_{i}}\delta_{J}(\tau_{i})}e'_{\tau_{i}},  \end{align*}and
  these two expressions are equal by (1). In the same fashion, we find
  that \begin{align*}N(f_{\mathcal{N}}(E_{\tau_i}))=f_{\mcal{N}}(N(E_{\tau_i}))=0, \end{align*}while \begin{align*}f_{\mcal{N}}(N(F_{\tau_i}))=f_{\mcal{N}}(0)=0,\end{align*}and \begin{align*}N(f_{\mcal{N}}(F_{\tau_i}))=N(u^{p\beta_{\tau_i}\delta_J(\tau_i)}f'_{\tau_i})=
  -\frac{(b)_{i-1}}{(a)_{i-1}}i_{\tau_{i-1}}\lambda_{\tau_{i-1}}u^{pi_{\tau_{i-1}}-p
    \alpha_{\tau_{i}}+p\beta_{\tau_i}\delta_J(\tau_i)}e'_{\tau_{i}}. \end{align*}If
  $\tau_i\notin J$, this expression is $0$ by (1). On the other hand,
  if $\tau_i\in J$, then the exponent of $u$ in this expression is
  $p(i_{\tau_{i-1}}-\alpha_{\tau_i}+\beta_{\tau_i})=pe$ by (6), so the
  expression is again $0$, as required.
\item Finally, that $f_{\mcal{M}}$ and $f_{\mcal{N}}$ commute with
  $\hat{g}$ follows directly from the definitions of $\alpha_{\tau_i}$
  and $\beta_{\tau_i}$.
\end{itemize}
It is clear from the construction that the maps $f_{\mathcal{M}''}$,
$f_{\mathcal{M}'}$, $f_{\mathcal{N}''}$ and $f_{\mathcal{N}'}$ are
nonzero. Since $T_{st}^*$ is faithful, the maps $T_{st}^{*}(f_{\mathcal{M}''})$,
$T_{st}^{*}(f_{\mathcal{M}'})$, $T_{st}^{*}(f_{\mathcal{N}''})$ and
$T_{st}^{*}(f_{\mathcal{N}'})$ are all nonzero, and are thus
isomorphisms (as they are maps between one-dimensional $E$-vector
spaces). The result follows.\end{proof}

\subsection{Weights and types}\label{reducibletypes}We recall some
definitions and results from \cite{dia05}. Fix, as ever, $\rho\sim
\bigl(\begin{smallmatrix}\psi_1&*\\0&\psi_2\end{smallmatrix}\bigr)$. We make
the following definitions:

\begin{defn}A weight $\sigma_{\vec{a},\vec{b}}$ is \emph{compatible} with $\rho$ (via $J$) if and only
if there exists a subset $J\in S$ so that
$$\psi_1|_{I_{K_0}}=\prod_{\tau\in S}\omega_\tau^{a_\tau}\prod_{\tau\in
J}\omega_\tau^{b_\tau},\ \psi_2|_{I_{K_0}}=\prod_{\tau\in
S}\omega_\tau^{a_\tau}\prod_{\tau\notin J}\omega_\tau^{b_\tau} $$
\end{defn}Suppose that these equations hold. We define
$$c_{\tau_i}=\left\{\begin{array}{ll}
                 {b_{\tau_i}}-\delta_J(\tau_{i+1})&\text{ if }\tau_i\in J \\
                 p-b_{\tau_i}-\delta_J(\tau_{i+1})&\text{ if }\tau_i\notin J \\
               \end{array}\right.$$where $\delta_J$ is the characteristic function of
               $J$. Define a character $\chi_{\vec{a},\vec{b},J}$ by
               $$\chi_{\vec{a},\vec{b},J}=\prod_{\tau_i\in
                 S}\omega_{\tau_i}^{a_{\tau_i}}\prod_{\tau_i\notin
                 J}\omega_{\tau_i}^{b_{\tau_i}-p}.$$Suppose that the
               $c_\tau$ are not all equal to either $0$ or $p-1$. Then
               we define a representation $I_{\vec{a},\vec{b},J}$ of
               $\GL_2(k)$ and a type $\tau_{\vec{a},\vec{b},J}$
               by $$I_{\vec{a},\vec{b},J}=I\left(\tilde{\chi}_{\vec{a},\vec{b},J},\tilde{\chi}_{\vec{a},\vec{b},J}\prod_{\tau\in
                   S}\tilde{\omega}_\tau^{c_\tau}\right)$$
$$\tau_{\vec{a},\vec{b},J}=\tilde{\chi}_{\vec{a},\vec{b},J}\oplus\tilde{\chi}_{\vec{a},\vec{b},J}\prod_{\tau\in
S}\tilde{\omega}_\tau^{c_\tau}$$Note that if $\rhobar$ is compatible
with $\sigma_{\vec{a},\vec{b}}$, then a lift of type $J$ is precisely
a lift of type $\tau_{\vec{a},\vec{b},J}$ with specified determinant.
\begin{prop}\label{typeswts}Suppose that $\sigma_{\vec{a},\vec{b}}$ is regular. If $\rho$ is compatible with $\sigma_{\vec{a},\vec{b}}$ via $J$, then $\rho$ is compatible with precisely one of the
Jordan-H\"{o}lder factors of the reduction mod $p$ of $I_{\vec{a},\vec{b},J}$,
and that factor is isomorphic to $\sigma_{\vec{a},\vec{b}}$.
\end{prop}
\begin{proof}We use the explicit computations of \cite{dia05}. Firstly, note
that reduction mod $p$ and the notion of compatibility both commute with
twisting, so we may replace $\rho$ by
$\rho\otimes\chi_{\vec{a},\vec{b},J}^{-1}$. By Proposition 1.1 of
\cite{dia05}, we have $\overline{I}_{\vec{a},\vec{b},J}\sim\oplus_{K\subset
S}\sigma_{\vec{a}_K,\vec{b}_K}$ where $a_K$ and $b_K$ are defined as follows:
$$a_{K,\tau_i}=\left\{\begin{array}{ll}
                 0&\text{ if }\tau_i\in K \\
                 c_{\tau_i}+\delta_K(\tau_{i+1})&\text{ if }\tau_i\notin K \\
               \end{array}\right.$$
$$b_{K,\tau_i}=\left\{\begin{array}{ll}
                 {c_{\tau_i}}+\delta_K(\tau_{i+1})&\text{ if }\tau_i\in K \\
                 p-c_{\tau_i}-\delta_K(\tau_{i+1})&\text{ if }\tau_i\notin K \\
               \end{array}\right.$$ By the definition of the $c_\tau$, we see at
               once that  $\sigma_{\vec{a}_J,\vec{b}_J}=\sigma_{\vec{a},\vec{b}}$, and in
               fact $$\psi_1|_{I_{K_0}}=\prod_{\tau\in S}\omega_\tau^{a_{J,\tau}}\prod_{\tau\in J}\omega_\tau^{b_{J,\tau}},\ \psi_2|_{I_{K_0}}=\prod_{\tau\in S}\omega_\tau^{a_{J,\tau}}\prod_{\tau\notin J}\omega_\tau^{b_{J,\tau}}.$$ 

If
$\rho$ is compatible with another Jordan-H\"{o}lder factor, there are  subsets $J'$,
$K'\subset S$, $J'\neq J$ such that
$$\psi_1|_{I_{K_0}}=\prod_{\tau\in S}\omega_\tau^{a_{J,\tau}}\prod_{\tau\in
J}\omega_\tau^{b_{J,\tau}}=\prod_{\tau\in
S}\omega_\tau^{a_{J',\tau}}\prod_{\tau\in K'}\omega_\tau^{b_{J',\tau}},$$
$$\psi_2|_{I_{K_0}}=\prod_{\tau\in S}\omega_\tau^{a_{J,\tau}}\prod_{\tau\notin
J}\omega_\tau^{b_{J,\tau}}=\prod_{\tau\in
S}\omega_\tau^{a_{J',\tau}}\prod_{\tau\notin K'}\omega_\tau^{b_{J',\tau}}.$$
Using the formulae above, the first equation simplifies to
$$\prod_{\tau_i\in S}\omega_{\tau_i}^{c_{\tau_i}+\delta_J(\tau_{i+1})}=\prod_{\tau_i\in (J'\cap K')\cup(J'^c\cap K'^c)}\omega_{\tau_i}^{c_{\tau_i}+\delta_{J'}(\tau_{i+1})}\prod_{\tau_{i+1}\in K'\cap J'^c}\omega_{\tau_i}.$$
By the assumption that $\sigma_{\vec{a},\vec{b}}$ is regular, we have $1\leq
c_{\tau_i}\leq p-2$ and $2\leq
c_{\tau_i}+\delta_{J}(\tau_{i+1})\leq p-2$ for each $i$. Then we see that we can equate the
exponents of $\omega_{\tau_i}$ on each side of each equation, and we easily
obtain $(J'\cap K')\cup(J'^c\cap K'^c)=S$, whence $J'=K'$. But then the
equation becomes
$$\prod_{\tau_i\in S}\omega_{\tau_i}^{\delta_J(\tau_{i+1})}=\prod_{\tau_i\in
S}\omega_{\tau_i}^{\delta_{J'}(\tau_{i+1})},$$whence $J=J'$, a
contradiction.
\end{proof}

\begin{rem}\label{weaklyregular}Note that it follows from the formulae in the proof of Proposition \ref{typeswts} that if $\sigma_{\vec{a},\vec{b}}$ is regular, then all the Jordan-H\"{o}lder factors of the reduction mod $p$ of $I_{\vec{a},\vec{b},J}$ are weakly regular.\end{rem}

\begin{prop}\label{jhcompat}Let $\theta_{1}$, $\theta_{2}$ be two tamely ramified characters of $I_{K_{0}}$ which extend to $G_{K_{0}}$. If $\rho$ has a potentially Barsotti-Tate lift (with determinant equal to a finite order character times the $p$-adic cyclotomic character) of type
$\theta_1\oplus\theta_2$, then $\rho$ is compatible with some weight
occurring in the mod $p$ reduction of $I(\theta_1,\theta_2)$.\end{prop}
\begin{proof}This follows easily from consideration of the possible
  Breuil modules corresponding to the $\pi_{L}$-torsion in the
  $p$-divisible group of such a lift (where the corresponding Galois
  representation is valued in $\bigO_{L}$, and $\pi_{L}$ is a
  uniformiser of $\bigO_{L}$). The case $\theta_1=\theta_2$ is easier,
  so from now on we assume that $\theta_1\neq\theta_2$. The
  $\pi_L$-torsion must contain a closed sub-group-scheme (with descent
  data) with generic fibre $\psi_1$. Suppose that this group scheme
  corresponds to a Breuil module with descent data $\M$. Then we can
  choose a basis so that $\M$ takes the following form:

$$\M^{\tau_i}=E[u]/u^{ep}\cdot x_{\tau_i} $$
$$\M_1^{\tau_i}=E[u]/u^{ep}\cdot u^{r_i}x_{\tau_i}$$
$$\phi_1(u^{r_i} x_{\tau_i})=(a^{-1})_ix_{\tau_{i+1}}$$
$$\hat{g}(x_{\tau_i})=\theta^i(g)x_{\tau_i}$$

Here $0\leq r_i\leq e$ is an integer, and $\theta^i:\Gal(K/K_0)\to
E^\times$ is a character. Now, by Corollary 5.2 of \cite{geesavquatalg},
because the lift is of type $\theta_1\oplus\theta_2$, we
must have $\theta^i=\theta_1$ or $\theta_2$ for each $i$ (here and below we denote
the reduction mod $p$ of the $\theta_{i}$ by the same symbol). Define
subsets $Y$, $Z$ by $$Y=\{\tau_i\in S|
\theta^i\neq\theta^{i+1}\},$$ $$Z=\{\tau_i\in S|
\theta^i=\theta_1\}.$$

Because $\theta_1\neq\theta_2$, if $i\in Y$ then the compatibility of the $\phi_1$- and $\Gal(K/K_{0})$-actions determines $r_i$ uniquely, and if $i\in Y^c$ then we can take either $r_i=0$ or $r_i=e$. Having written down all possible $\M$, we now need to determine their generic fibres. This is a straightforward calculation using Example 3.7 of \cite{sav06}. Without loss of generality, we may twist and assume that $\theta_1=\prod_{\tau_i\in S}\omega_{\tau_i}^{c_i}$, $\theta_2=1$, with $0\leq c_i\leq p-1$. Then one easily obtains $$\psi_1|_{I_{K_0}}=\omega_{\tau_1}^{m_1+n_1}\prod_{\tau_i\in \{Y^c|r_i=e\}}\omega_{\tau_i}\prod_{\tau_i\in Y\cap Z}\omega_{\tau_i},$$where
$$m_1=\left\{\begin{array}{ll}
                 0&\text{ if }\tau_1\notin Z \\
                 c_1+pc_r+\dots+p^{r-1}c_2&\text{ if }\tau_1\in Z \\
               \end{array}\right.$$
\begin{align*}n_1=&\frac{1}{e}\sum_{\tau_i\in Y\cap
  Z^c}p^{r-i}(p^ic_1+p^{i+1}c_r+\dots+p^rc_i+c_{i+1}+\dots+p^{i-1}c_2)\\
&-\frac{1}{e}\sum_{\tau_i\in Y\cap
  Z}p^{r-i}(p^ic_1+p^{i+1}c_r+\dots+p^rc_i+c_{i+1}+\dots+p^{i-1}c_2).\end{align*}Now,
consider the coefficient of $c_1$ in $n_1$. The sets $Y\cap Z^c$ and
$Y\cap Z$ have equal cardinality, so this coefficient is in fact
zero. Thus the coefficient of $c_1$ in $m_1+n_1$ is $1$ if $\tau_1\in
Z$, and $0$ otherwise. By cyclic symmetry, we
obtain $$\psi_1|_{I_{K_0}}=\prod_{\tau_i\in
  Z}\omega_{\tau_i}^{c_i}\prod_{\tau_{i}\in X}\omega_{\tau_i},$$
where $$X=\{\tau_i\in Y^c|r_i=e\}\cup(Y\cap Z).$$

We wish to show that $\rho$ is compatible with some weight in the reduction mod $p$ of $I(\theta_1,\theta_2)$. It is easy to check that the determinant of $\rho$ is correct, so it suffices to examine $\psi_1$; in the notation of Proposition \ref{typeswts}, we see that
$\rho$ is compatible with $\sigma_{\vec{a}_K,\vec{b}_K}$ via $L$ if and only if $$\psi_1|_{I_{K_0}}=\prod_{\tau_i\in (K^{c}\cap L)\cup(K\cap L^c)}\omega_{\tau_i}^{c_i+\delta_{K^{c}}(\tau_{i+1})}\prod_{\tau_i\in S}\omega_{\tau_i}^{\delta_{K\cap L}(\tau_{i+1})}$$(note that our convention that $\theta_{2}=1$ causes $K^{c}$ to appear in this formula rather than $K$).

The result now follows upon taking, for
example, $$K=\{\tau_i|\tau_{i-1}\in(X^c\cap Y^c\cap Z)\cup(X\cap
Y^c\cap Z^c)\}$$
and $$L=(K^c\cap Z)\cup(K\cap Z^c).$$\end{proof}

\begin{prop}\label{prop:liftimpliesmodel}
	Suppose that $\sigma_{\vec{a},\vec{b}}$ is regular. If $\rho$ is compatible with $\sigma_{\vec{a},\vec{b}}$ via $J$, and $\rho$ has a lift of type $J$, then $\rho$ has a model of type $J$.
\end{prop}
\begin{proof}\label{pf:liftimpliesmodel}This follows from similar considerations to those involved in the proof of Proposition \ref{jhcompat}. Consider the $\pi_{L}$-torsion in the $p$-divisible group corresponding to the lift of type $J$. It contains a closed sub-group-scheme (with descent data) with generic fibre $\psi_1$. Suppose that this group scheme corresponds to a Breuil module with descent data $\M$. Then we can choose a basis so that $\M$ takes the following form:

$$\M^{\tau_i}=E[u]/u^{ep}\cdot x_{\tau_i} $$
$$\M_1^{\tau_i}=E[u]/u^{ep}\cdot u^{r_i}x_{\tau_i}$$
$$\phi_1(u^{r_i} x_{\tau_i})=(a^{-1})_ix_{\tau_{i+1}}$$
$$\hat{g}(x_{\tau_i})=\theta^i(g)x_{\tau_i}$$
Again, by Corollary 5.2 of \cite{geesavquatalg} and the definition of a lift of type $J$, for each $i$ we must have $\theta^i=\theta_1$ or $\theta^{i}=\theta_{2}$ where $$\theta_{1}=\prod_{\tau\in S}\omega_{\tau}^{a_{\tau}}\prod_{\tau\in J}\omega_{\tau}^{b_{\tau}-p},$$  $$\theta_{2}=\prod_{\tau\in S}\omega_{\tau}^{a_{\tau}}\prod_{\tau\in J^{c}}\omega_{\tau}^{b_{\tau}-p}.$$ Note that $\psi_1|_{I_{K_0}}=\theta_{1}\prod_{\tau_{i}\in S}\omega_{\tau_{i}}^{\delta_{J}(\tau_{i+1})}$. Without loss of generality, we can twist so that $\theta_1=\prod_{\tau_i\in S}\omega_{\tau_i}^{c_i}$, $\theta_2=1$, with $0\leq c_i\leq p-1$. Then we obtain $$\theta_{1}=\theta_{1}\theta_{2}^{-1}=\prod_{\tau_{i}\in J}\omega_{\tau_{i}}^{b_{\tau_{j}}-\delta_{J}(\tau_{i+1})}\prod_{\tau_{i}\in J^{c}}\omega_{\tau_{i}}^{p-b_{\tau_{i}}-\delta_{J}(\tau_{i+1})}. $$Since $0\leq c_{i}\leq p-1$ and $2\leq b_{\tau_{i}}\leq p-2$, we obtain $$c_{i}=\left\{\begin{array}{ll}
                 {b_{\tau_i}}-\delta_J(\tau_{i+1})&\text{ if }\tau_i\in J \\
                 p-b_{\tau_i}-\delta_J(\tau_{i+1})&\text{ if }\tau_i\notin J \\
               \end{array}\right.$$Note that this implies that $2\leq c_{i}+\delta_{J}(\tau_{i+1})\leq p-2$. Now, using the same definitions of $X$, $Y$ and $Z$ as in the proof of Proposition \ref{jhcompat}, we can compare the two expressions we have for $\psi_{1}|_{I_{K_0}}$ to obtain $$\prod_{\tau_{i}\in S}\omega_{\tau_{i}}^{c_{i}+\delta_{J}(\tau_{i+1})}=\prod_{\tau_i\in Z}\omega_{\tau_i}^{c_i}\prod_{\tau_{i}\in X}\omega_{\tau_i}. $$ Since $2\leq c_{i}+\delta_{J}(\tau_{i+1})\leq p-2$, this gives $Z=S$, and $X=\{\tau_{i}|\tau_{i+1}\in J\}$. Since $Z=S$, we have $Y=\emptyset$, and thus the fact that  $X=\{\tau_{i}|\tau_{i+1}\in J\}$ means that $\M$ is in fact of class $J$. It is then clear that the $\pi_{L}$-torsion is a model of $\rho$ of type $J$, as required.
	
\end{proof}

\section{Local analysis - the irreducible case}\label{irreducible}\subsection{}We now
prove the analogues of some of the results of section \ref{reducible} in the case
where $\rho$ is irreducible. 

We assume that $\rho$ is irreducible from now on. In addition to the
assumptions made at the beginning of section \ref{reducible}, we now also assume that
$\F_{p^2}\subset E$, where $\rho:G_{K_0}\to\GL_2(E)$. Let $k'$ be the
(unique) quadratic extension of $k$.

Label the embeddings $k'\into\Fpbar$ as $S'=\{\sigma_{i}\}$, $0\leq i\leq 2r-1$, so that $\sigma_{i+1}=\sigma_{i}\circ\phi^{-1}$, and $\sigma_{i}|_{k}=\tau_{\pi(i)}$, where $\pi:\Z/2r\Z\to\Z/r\Z$ is the natural surjection. For simplicity of notation we will sometimes refer to the elements of $S'$ as elements of $\Z/2r\Z$, and the elements of $S$ as elements of $\Z/r\Z$.

Recall that we say that a subset $H\subset S'$ is a \emph{full subset} if $|H|=|\pi(H)|=r$.

\begin{defn}We say that $\rho$ is \emph{compatible} with a weight $\sigma_{\vec{a},\vec{b}}$ (via $J$)
if there exists a full subset $J\subset S'$ such that
$$\rho|_{I_{K'_0}}\sim\prod_{\sigma\in S'}\omega_\sigma^{a_\sigma}\left(\begin{array}{cc}
  \prod_{\sigma\in J}\omega_\sigma^{b_\sigma}  & 0 \\
  0 & \prod_{\sigma\notin J}\omega_\sigma^{b_\sigma} \\
\end{array}\right),$$where we write $a_{\sigma}$, $b_{\sigma}$ for $a_{\pi(\sigma)}$, $b_{\pi(\sigma)}$ respectively.\end{defn} Note that the predicted set of weights $W(\rhobar)$ is
just the set of compatible weights; this is one way in which the irreducible
case is simpler than the reducible one.

 Given a regular weight $\sigma_{\vec{a},\vec{b}}$ and a full subset $J\subset S'$, we wish to define a representation and a type. Let $K_{J}=\pi(J\cap\{1,\dots,r\})$. Then let 

$$c_{i}=\left\{\begin{array}{llll}
                 {b_{i}}+\delta_{K_{J}}({1})-1&\text{ if }0=i\in K_{J} \\
                 p-b_{i}+\delta_{K_{J}}({1})-1&\text{ if }0=i\notin K_{J} \\
				{b_{i}}-\delta_{K_{J}}({i+1})&\text{ if }0\neq i\in K_{J} \\
                 p-b_{i}-\delta_{K_{J}}({i+1})&\text{ if }0\neq i\notin K_{J} \\
               \end{array}\right.$$
Define a character $$\psi_{\vec{a},\vec{b},J}=\omega_{\tau_{0}}^{-\delta_{K_{J}}(1)}\prod_{\tau\in S}\omega_{\tau}^{a_{\tau}}\prod_{\tau\notin {K}_{J}}\omega_{\tau}^{b_{\tau}-p}.$$Then we define

$$I'_{\vec{a},\vec{b},J}=\Theta\left(\tilde{\psi}_{\vec{a},\vec{b},J}\tilde{\omega}_{\sigma_{r}}\prod_{i=1}^{r}\tilde{\omega}_{\sigma_{i}}^{c_{i}}\right)$$

$$\tau'_{\vec{a},\vec{b},J}=\tilde{\psi}_{\vec{a},\vec{b},J}\tilde{\omega}_{\sigma_{r}}\prod_{i=1}^{r}\tilde{\omega}_{\sigma_{i}}^{c_{i}}\oplus\left(\tilde{\psi}_{\vec{a},\vec{b},J}\tilde{\omega}_{\sigma_{r}}\prod_{i=1}^{r}\tilde{\omega}_{\sigma_{i}}^{c_{i}}\right)^{p^r}$$

\begin{prop}\label{typeswts2}Recall that $\sigma_{\vec{a},\vec{b}}$ is regular. If $\rho$ is compatible with $\sigma_{\vec{a},\vec{b}}$ via $J$, then $\rho$ is compatible with precisely one of the
Jordan-H\"{o}lder factors of the reduction mod $p$ of $I'_{\vec{a},\vec{b},J}$,
and that factor is isomorphic to $\sigma_{\vec{a},\vec{b}}$.
\end{prop}
\begin{proof}We may twist and assume without loss of generality that
  $\psi_{\vec{a},\vec{b},J}=1$. Then by Proposition 1.3 of
  \cite{dia05} (note here that Diamond's conventions on the numbering
  of the elements of $S'$ are the opposite of ours, so that his
  $\sigma_{i}$ is our $\sigma_{-i}$), the Jordan-H\"{o}lder factors of
  the reduction mod $p$ of $I'_{\vec{a},\vec{b},J}$ are
  $\{\sigma_{\vec{a}_{K},\vec{b}_{K}}\}_{{K\subset S}}$, where
	$$a_{K,\tau_{i}}=\left\{\begin{array}{llll}
	                 \delta_{K}({1})&\text{ if }0=i\in K\\
	                c_{i}+1&\text{ if }0=i\notin K \\
				0&\text{ if }0\neq i\in K \\
	                c_{i}+\delta_{K}({i+1})&\text{ if }0\neq i\notin K \\
	               \end{array}\right.$$	
				$$b_{K,\tau_{i}}=\left\{\begin{array}{llll}
				                 {c_{i}}+1-\delta_{K}({1})&\text{ if }0=i\in K \\
				                 p-c_{i}+\delta_{K}({1})-1&\text{ if }0=i\notin K \\
								{c_{i}}+\delta_{K}({i+1})&\text{ if }0\neq i\in K \\
				                 p-c_{i}-\delta_{K}({i+1})&\text{ if }0\neq i\notin K \\
               \end{array}\right.$$From the definition of the $c_{i}$ and of $\psi_{\vec{a},\vec{b},J}$, we have $\sigma_{\vec{a}_{K_{J}},\vec{b}_{K_{J}}}=\sigma_{\vec{a},\vec{b}}$. Suppose that $\rho$ is compatible with $\sigma_{\vec{a}_{K'},\vec{b}_{K'}}$ via $J'$. Then, replacing $J'$ by $(J')^{c}$ if necessary, we must have $$\prod_{i\in S'}\omega_{\sigma_{i}}^{a_{K_{J},i}}\prod_{i\in J}\omega_{\sigma_{i}}^{b_{K_{J},i}}=\prod_{i\in S'}\omega_{\sigma_{i}}^{a_{K',i}}\prod_{i\in J'}\omega_{\sigma_{i}}^{b_{K',i}}.	$$ Using the formulae above, this becomes \begin{align}\label{407eqn}
              \omega_{\sigma_{0}}^{\delta_{J',K'}(1)}\omega_{\sigma_{r}}^{\delta_{J',K'}(r+1)}\prod_{i\in T'}\omega_{\sigma_{i}}^{c_{i}+\delta_{K'}(i+1)}\prod_{i\in S'}\omega_{\sigma_{i}}^{\delta_{{J'\cap\pi^{-1}((K')^{c})}}(i+1)}\notag \\	= \omega_{\sigma_{0}}^{\delta_{J,K_{J}}(1)}\omega_{\sigma_{r}}^{\delta_{J,K_{J}}(r+1)}\prod_{i\in T}\omega_{\sigma_{i}}^{c_{i}+\delta_{K_{J}}(i+1)}\prod_{i\in S'}\omega_{\sigma_{i}}^{\delta_{{J\cap\pi^{-1} (K_{J}^{c})}}(i+1)},
               \end{align}	where $$T=(J\cap\pi^{-1}(K_{J}))\cup(J^{c}\cap\pi^{-1}(K_{J}^{c}))=\{1,\dots,r\}, $$ $$T'=(J'\cap\pi^{-1}(K'))\cup((J')^{c}\cap\pi^{-1}((K')^{c})), $$ $$\delta_{J,K_{J}}(i+1)=\left\{\begin{array}{ll}
			                 1-\delta_{K_{J}}({i+1})&\text{ if }i\in T\\
			              \delta_{K_{J}}(i+1)&\text{ if }i\notin T, \\ \end{array}\right. $$ 	$$\delta_{J',K'}(i+1)=\left\{\begin{array}{ll}
							                 1-\delta_{K'}({i+1})&\text{ if }i\in T'\\
							              \delta_{K'}(i+1)&\text{ if }i\notin T'. \\ \end{array}\right. $$Note that (since $\sigma_{\vec{a},\vec{b}}$ is regular) all the exponents on the right hand side of (\ref{407eqn}) are in the range $[0,p-1]$. On the left hand side, this is true except possibly for the exponents of $\omega_{\sigma_{0}}$, $\omega_{\sigma_{r}}$. Since $T=\{1,\dots,r\}$, it is easy to see that the only opportunity for this not to hold is for the exponent of $\omega_{\sigma_{0}}$ to be $p$ on the left hand side and $0$ on the right hand side. However, in order for the exponent of  $\omega_{\sigma_{0}}$ to be $p$ on the left hand side we require $c_{0}=p-2$, which requires that $1\in K_{J}$. But then the exponent of $\omega_{\sigma_{0}}$ on the right hand side is 1, a contradiction.
							
							Thus we may equate exponents on each side of (\ref{407eqn}). In particular, if $i\neq 0$, we have (again because $\sigma_{\vec{a},\vec{b}}$ is regular) $c_{i}+\delta_{K_{J}}(i+1)\in [2,p-2]$, so that we must have $\{1,\dots,r-1\}\subset T'$. We also have $c_{0}\in[1,p-2]$. If $0\in T'$, we see that the exponent of $\omega_{\sigma_{0}}$ on the left hand side of (\ref{407eqn}) is $c_{0}+1+\delta_{J'\cap\pi^{-1}((K')^{c})}(1)=c_{0}+1$ (because $1\in T'$), which is at least 2. However the exponent of $\omega_{\sigma_{0}}$ on the right hand side of (\ref{407eqn}) is 0 or 1, as $0\notin T$, which is a contradiction. Thus $T'=T=\{1,\dots,r\}$.
							
							 Then (\ref{407eqn}) simplifies to $$\prod_{i=0}^{r-1}\omega_{\sigma_{i}}^{\delta_{K'}(i+1)}\prod_{i=r}^{2r-1}\omega_{\sigma_{i}}^{\delta_{(K')^{c}}(i+1)}=\prod_{i=0}^{r-1}\omega_{\sigma_{i}}^{\delta_{K_{J}}(i+1)}\prod_{i=r}^{2r-1}\omega_{\sigma_{i}}^{\delta_{K_{J}^{c}}(i+1)}, $$whence $K'=K_{J}$, as required. \end{proof}

\begin{rem}\label{weaklyregular2}Note that it follows easily from the formulae in the proof of Proposition \ref{typeswts2} that if $\sigma_{\vec{a},\vec{b}}$ is regular, then all the Jordan-H\"{o}lder factors of the reduction mod $p$ of $I'_{\vec{a},\vec{b},J}$ are weakly regular.\end{rem}

\begin{thm}\label{localdef2}Assume that $\sigma_{\vec{a},\vec{b}}$ is
  regular and that $\rho$ is compatible with
  $\sigma_{\vec{a},\vec{b}}$ via $J$. Then $\rho$ has a lift of type
  $\tau'_{\vec{a},\vec{b},J}$ which is not potentially ordinary. 
\end{thm}
\begin{proof}A simple computation shows that we in fact have $$\tau'_{\vec{a},\vec{b},J}=\prod_{\tau\in S}\omega_{\tau}^{a_{\tau}}\prod_{\sigma\in J}\omega_{\sigma}^{b_{\sigma}-p}\oplus\prod_{\tau\in S}\omega_{\tau}^{a_{\tau}}\prod_{\sigma\notin J}\omega_{\sigma}^{b_{\sigma}-p}. $$This means that we only need to make a very minor modification to the proof of Theorem \ref{localdef}. Let $K_{0}'=W(k')[1/p]$. Fix
$\pi'=(-p)^{1/(p^{2r}-1)}$, and let $K'=K_0'(\pi')$. Let $g_{\phi}$ be the nontrivial element of $\Gal(K'/K_{0})$ which fixes $\pi'$. It is clear from the proof of Theorem \ref{localdef} that for some choice of $a\in W(E)^{\times}$ the following object of $W(E)-\Mod^{1}_{cris,dd,K_{0}}$ provides us with the required lift.

	$$\M_J^{\sigma_i}=S_{K}\cdot e_{\sigma_i}+S_{K}\cdot f_{\sigma_i}$$
	$$\hat{g}_{\phi}(e_{\sigma_i})=f_{\sigma_{i+r}}$$ $$\hat{g}_{\phi}(f_{\sigma_i})=e_{\sigma_{i+r}}$$
	If $g\in\Gal(K'/K_{0}')$,
	 $$\hat{g}(e_{\sigma_i})=\left(\left(\prod_{\tau\in S}\widetilde{\omega}_{\tau}^{a_{\tau}}\prod_{\sigma\in J}\widetilde{\omega}_{\sigma}^{b_{\sigma}-p}\right)(g)\right)e_{\sigma_i}$$
	$$\hat{g}(f_{\sigma_i})=\left(\left(\prod_{\tau\in S}\widetilde{\omega}_{\tau}^{a_{\tau}}\prod_{\sigma\notin J}\widetilde{\omega}_{\sigma}^{b_{\sigma}-p}\right)(g)\right)f_{\sigma_i}$$
		If  $\sigma_{i+1}\in J$,
	$$\Fil^1\M_J^{\sigma_i}=\Fil^1S_{K}\cdot\M^{\sigma_i}_J+S_{K}\cdot f_{\sigma_i}$$ $$\phi(e_{\sigma_i})=(a^{-1})_i e_{\sigma_{i+1}}$$ $$\phi(f_{\sigma_i})=({a}^{-1})'_i pf_{\sigma_{i+1}}$$
		If  $\sigma_{i+1}\notin J$,
	$$\Fil^1\M_J^{\sigma_i}=\Fil^1S_{K}\cdot\M^{\sigma_i}_J+S_{K}\cdot e_{\sigma_i}$$ $$\phi(e_{\sigma_i})=({a}^{-1})_i p e_{\sigma_{i+1}}$$ $$\phi(f_{\sigma_i})=({a}^{-1})'_i f_{\sigma_{i+1}}$$Here the notation $(x)'_{i}$ means $x$ if $i=r+1$ and $1$ otherwise.

\end{proof}
\section{Global Results}\label{global}\subsection{}We now show how the local results obtained in the
previous sections can be combined with lifting theorems to prove results
about the possible weights of mod $p$ Hilbert modular forms. Firstly, we show that if $\rhobar$ is modular of some regular weight, then $\rhobar$ is compatible with that weight, by making use of Lemma \ref{lem:local langlands version of lifting} and Proposition \ref{jhcompat}. We then turn this analysis around. We take a conjectural regular weight $\sigma$ for $\rhobar$, and using modularity lifting theorems we produce a modular lift of $\rhobar$ of a specific type, which is enough to prove that $\rhobar$ is modular of weight $\sigma$ by Propositions \ref{typeswts} and \ref{typeswts2}.

Assume now that $F$ is a totally real field in which $p>2$ is unramified,
and that $\overline{\rho}:G_F\to\GL_2(E)$ is a continuous representation,
known to be modular, where $E$ is a finite extension of $\F_p$.

Let $W(\overline{\rho})$ be the conjectural set of Serre weights for
$\overline{\rho}$, as defined in Section \ref{2}. Recall that the elements of
$W(\overline{\rho})$ are just the tensor products of elements of
$W_v(\overline{\rho})$, for $v|p$, and that such elements are of the form
$\sigma_{\vec{a},\vec{b}}$, as described above. We say that a weight is (weakly) regular
if and only if it is a tensor product of (weakly) regular weights. 

The following argument is based on an argument of Michael Schein (c.f. Proposition 5.11 of \cite{sch06}), and is due to him in the case that $\overline{\rho}|_{G_{F_v}}$ is irreducible.

\begin{lemma}\label{compat} Suppose that $p\geq 3$, that $\overline{\rho}$ is modular of weight $\sigma=\otimes_v \sigma^v_{\vec{a},\vec{b}}$, and that $\sigma$ is weakly regular. Then for each $v$, either $\overline{\rho}|_{G_{F_v}}$ is compatible with $\sigma^v_{\vec{a},\vec{b}}$, or $\sigma^v_{\vec{a},\vec{b}}$ is not regular and $\overline{\rho}|_{G_{F_v}}$ is not compatible with any regular weight.\end{lemma}
\begin{proof}

Suppose firstly that $\overline{\rho}|_{G_{F_v}}$ is reducible. We
will assume for the rest of this proof that $F_v\neq\Q_p$; the
argument needed when $F_v=\Qp$ is slightly different, although much
simpler, and the result follows from Lemma 4.4.6 of \cite{gee061}. We
will also assume that there is at least one $b_{\tau_i}\ne 1$; the
case where all $b_{\tau_i}=1$ is much easier, and we leave it to the reader. Then for any type $\tau=\chi_1\oplus\chi_2$ (with $\chi_{1}\ne\chi_{2}$ tame characters of $I_{F_{v}}$ which extend to $G_{F_{v}}$) such that
$\sigma^v_{\vec{a},\vec{b}}$ occurs in the reduction of $I(\chi_1,\chi_2)$, it follows from Lemma \ref{lem:local langlands version of lifting} and Proposition \ref{jhcompat} that there must be a weight
$\sigma^v_{\vec{a'},\vec{b'}}$ in the reduction of $I(\chi_1,\chi_2)$ which is
compatible with $\overline{\rho}|_{G_{F_v}}$. Since we are
working purely locally, we return to the notation of section
\ref{reducibletypes}.

Twisting, we may without loss of generality suppose that $a_\tau=0$ for all
$\tau$. By Proposition 1.1 of \cite{dia05} (and the fact that $\sigma$
is weakly regular, with at least one $b_{\tau_i}\ne 1$) there is for each $J\subset S$ a
unique pair of characters $\prod_{\tau\in S}\tilde{\omega}^{c^J_\tau}$, $\prod_{\tau\in S}\tilde{\omega}^{d^J_\tau}$ (with $0\leq c^{J}_{\tau}, d^{J}_{\tau}\leq p-1$) such that if we define
$$\sigma^J=I\left(1,\prod_{\tau\in
S}\tilde{\omega}^{d^J_\tau}\right)\otimes\prod_{\tau\in S}\tilde{\omega}_{\tau}^{c^{J}_{\tau}}\circ\det$$then, with the same
notation for reductions as in \cite{dia05}, extended to be compatible with twisting, $\sigma^J_J\sim \sigma_{\vec{a},\vec{b}}$.
Then there must (by the argument above) be some subset $K_J\subset S$, such that
$\sigma^J_{K_J}$ is compatible with $\rho$. If $\sigma^J_{K_J}\sim
\sigma_{\vec{m}^J_{K_J},\vec{n}^J_{K_J}}$ this means that there must be a subset
$L_J\subset S$ such that $$\psi_1|_{I_{K_0}}=\prod_{\tau\in
S}\omega_\tau^{{m_{K_{J},\tau}^J}}\prod_{\tau\in L_J}\omega_\tau^{{n_{K_{J},\tau}^J}}.$$

By Proposition 1.1 of \cite{dia05}, this is equal to
$$\prod_{\tau_i\in S}\omega_{\tau_i}^{c^J_{\tau_i}}\prod_{\tau_i\in
L_J\cap K_J^c}\omega_{\tau_i}^p\prod_{\tau_i\in(L_J\cap K_J)\cup(L_J^c\cap
K_J^c)}\omega_{\tau_i}^{d^J_{\tau_i}+\delta_{K_J}(\tau_{i+1})}.$$

Now, since
$\sigma^J_J\sim \sigma_{\vec{a},\vec{b}}$, we have $$\prod_{\tau_i\in
S}\omega_{\tau_i}^{c^J_{\tau_i}}\prod_{\tau_i\notin
J}\omega_{\tau_i}^{d^J_{\tau_i}+\delta_J(\tau_{i+1})}=\prod_{\tau_i\in
S}\omega_{\tau_i}^{a_{\tau_i}}=1,$$

by the assumption that $a_{\tau}=0$ for all $\tau$, so that in fact $$\psi_1|_{I_{K_0}}=\prod_{\tau_i\in
J^c}\omega_{\tau_i}^{-(d^J_{\tau_i}+\delta_J(\tau_{i+1}))}\prod_{\tau_i\in
L_J\cap K_J^c}\omega_{\tau_i}^p\prod_{\tau_i\in(L_J\cap K_J)\cup(L_J^c\cap
K_J^c)}\omega_{\tau_i}^{d^J_{\tau_i}+\delta_{K_J}(\tau_{i+1})}.$$

Since $\sigma^J_J\sim \sigma_{\vec{a},\vec{b}}$, we have
$$d^J_{\tau_i}=\left\{\begin{array}{ll}
                 b_{\tau_i}-\delta_J(\tau_{i+1})&\text{ if }\tau_i\in J \\
                 p-b_{\tau_i}-\delta_J(\tau_{i+1})&\text{ if }\tau_{i}\notin J \\
               \end{array}\right.$$
Substituting, we see that $$\psi_1|_{I_{K_0}}=\prod_{\tau_i\in
(T_J\cap J)\cup(T_J^c\cap J^c)}\omega_{\tau_i}^{b_{\tau_i}}\prod_{\tau_i\in
S}\omega_{\tau_i}^{\delta_{L_J\cap K_J^c}(\tau_{i+1})-\delta_{T^c_J\cap
J^c}(\tau_{i+1})} \prod_{\tau_i\in
T_J}\omega_{\tau_i}^{\delta_{K_J}(\tau_{i+1})-\delta_{J}(\tau_{i+1})},$$ where we write $T_J=(K_J\cap L_J)\cup(K^c_J\cap L^c_J)$.

Putting $J=S$, we obtain \begin{align}\psi_1|_{I_{K_0}}&=\prod_{\tau_i\in
T_S}\omega_{\tau_i}^{b_{\tau_i}}\prod_{\tau_i\in
S}\omega_{\tau_i}^{\delta_{L_S\cap K_S^c}(\tau_{i+1})} \prod_{\tau_i\in
T_S}\omega_{\tau_i}^{\delta_{K_S}(\tau_{i+1})-1}\notag\\ &=\prod_{\tau_{i}\in T_{S}}\omega_{\tau_{i}}^{b_{\tau_{i}}-\delta_{K_{S}^{c}\cap L_{S}^{c}}(\tau_{i+1})}\prod_{\tau_{i}\in T^{c}_{S}}\omega_{\tau_{i}}^{\delta_{L_{S}\cap K_{S}^{c}}(\tau_{i+1})}.\label{Sexpression} \end{align}

Now, suppose that $\sigma_{\vec{a},\vec{b}}$ is \emph{not} compatible
with $\rho$, so that for all $J$ we have $K_{J}\neq J$.  We can
uniquely write $$\psi_1|_{I_{K_0}}=\prod_{\tau_i\in
  S}\omega_{\tau_i}^{c_{\tau_i}}$$ with $0\leq c_{\tau_i}\leq p-1$ not
all equal to $p-1$ (in fact, an examination of the product just
written shows that the exponents are already in this range). Examining
the formula just established, we see that after possibly exchanging
$\psi_1$ and $\psi_2$ (which we can do, as the definition of
``compatible'' is unchanged by this exchange), there must be some $j$
such that $b_{\tau_j}\neq 1$, $c_{\tau_j}=b_{\tau_j}-1$, $\tau_j\in
T_S$, and $\tau_{j+1}\in K_S^c\cap L_S^c\subset T_S$ (else
$\rho$ would be compatible with $\sigma_{\vec{a},\vec{b}}$).

Now take $J=\{\tau_{j}\}$, so that \begin{align}\psi_1|_{I_{K_0}}&=\prod_{\tau_i\in
(T_{\{\tau_{j}\}}\cap \{\tau_{j}\})\cup(T_{\{\tau_{j}\}}^{c}\cap \{\tau_{j}\}^c)}\omega_{\tau_i}^{b_{\tau_i}}\prod_{\tau_i\in
S}\omega_{\tau_i}^{\delta_{L_{\{\tau_{j}\}}\cap K_{\{\tau_{j}\}}^{c}}(\tau_{i+1})-\delta_{T^c_{\{\tau_{j}\}}\cap
\{\tau_{j}\}^c}(\tau_{i+1})}\cdot\notag\\& \prod_{\tau_i\in
T_{\{\tau_{j}\}}}\omega_{\tau_i}^{\delta_{K_{\{\tau_{j}\}}}(\tau_{i+1})-\delta_{\{\tau_{j}\}}(\tau_{i+1})}.\label{jexpression}\end{align}It is easy to see that the exponent of $\omega_{\tau_{i}}$ in this product is always between 0 and $p-1$, unless $i=j-1$ or $i=j$. If the exponent is always between $0$ and $p-1$, then we have a contradiction, because we already know that $c_{\tau_j}=b_{\tau_j}-1$, but from (\ref{jexpression}) we see that the exponent of $\omega_{\tau_{j}}$ can only be $0$, $b_{\tau_{j}}$ or $b_{\tau_{j}}+1$.

So, at least one of the exponents of $\omega_{\tau_{j-1}}$ and
$\omega_{\tau_{j}}$ must be $-1$ or $p$. We now analyse when this can
occur. It's easy to see that the exponent of $\omega_{\tau_{j}}$ is
$-1$ if and only if $\tau_{j}\notin T_{\{\tau_{j}\}}$ and
$\tau_{j+1}\in L^{c}_{\{\tau_{j}\}}\cap K_{\{\tau_{j}\}}$, and it is
$p$ if and only if $b_{\tau_{j}}=p-1$, $\tau_{j}\in T_{\{\tau_{j}\}}$
and $\tau_{j+1}\in L_{\{\tau_{j}\}}\cap K_{\{\tau_{j}\}}$. Similarly,
the exponent of $\omega_{\tau_{j-1}}$ is $-1$ if and only if
$\tau_{j-1}\in T_{\{\tau_{j}\}}$ and $\tau_{j}\in
L^{c}_{\{\tau_{j}\}}\cap K^{c}_{\{\tau_{j}\}}$, and it is $p$ if and
only if $b_{\tau_{j-1}}=p-1$, $\tau_{j-1}\in T^{c}_{\{\tau_{j}\}}$ and
$\tau_{j}\in L_{\{\tau_{j}\}}\cap K^{c}_{\{\tau_{j}\}}$. Thus it is
impossible for both exponents to be $p$, or both to be $-1$.

Suppose now that the exponent of $\omega_{\tau_{j}}$ in (\ref{jexpression}) is $-1$. If we multiply each of the expressions (\ref{Sexpression}), (\ref{jexpression}) by $\omega_{\tau_{j}}$, write each side as a product $\prod_{\tau}\omega_{\tau}^{n_{\tau}}$ with $0\leq n_{\tau}\leq p-1$ and equate coefficients of $\omega_{\tau_{j}}$ in the resulting expression, we obtain $b_{\tau_{j}}=0$ or $1$ (the second case only a possibility when the exponent of $\omega_{\tau_{j-1}}$ in (\ref{jexpression}) is $p$), a contradiction.

Suppose that the exponent of $\omega_{\tau_{j}}$ in
(\ref{jexpression}) is $p$. Then we again easily see that
$p-2=b_{\tau_{j}}-1=0$ or $1$. Thus $p-2=1$, and we additionally need
to have
$(T_{\{\tau_j\}}\cap\{\tau_j\}\cup(T_{\{\tau_j\}}^c\cap\{\tau_j\}^c=S$, so
that $T_{\{\tau_j\}}=\{\tau_j\}$. But for the exponent of
$\omega_{\tau_j}$ to be $p$ we need that $\tau_{j+1}\in
L_{\{\tau_j}\}\cap K_{\{\tau_j\}}\subset T_{\{\tau_j\}}$, a contradiction.

Suppose that the exponent of $\omega_{\tau_{j-1}}$ in
(\ref{jexpression}) is $p$. Then in the same fashion we obtain
$b_{\tau_{j}}-1$=$0$, or $1$. The only possibility is that $b_{\tau_{j}}=2$, when we in addition (in order that the necessary carrying should occur) require that $b_{\tau_{i}}=p-1$ for all $i\neq j$.

Finally, suppose that the exponent of $\omega_{\tau_{j-1}}$ in (\ref{jexpression}) is $-1$. Multiply each of  (\ref{Sexpression}), (\ref{jexpression}) by $\omega_{\tau_{j-1}}$. Then we see that the only way for equality to hold is again if $b_{\tau_{i}}=p-1$ for all $i \neq j$.

So, we have deduced that $b_{\tau_{i}}=p-1$ for all $i \neq j$, so that $\sigma_{\vec{a},\vec{b}}$ is certainly not regular. It now remains to show that  $\rho$ is not compatible with any regular weight. Examining the above argument, we see that we have in fact deduced that (again, after possibly exchanging $\psi_{1}$, $\psi_{2}$) $$\psi_{1}|_{I_{K}}=\omega_{\tau_{j}}^{b_{\tau_{j}}-1}\prod_{i\neq j}\omega_{\tau_{i}}^{p-1}, $$ $$\psi_{2}|_{I_{K}}=\omega_{\tau_{j}}, $$ with $2\leq b_{\tau_{j}}\leq p-1$.

If $\rho$ is compatible with some regular weight, then we have by definition that $$\psi_{1}|_{I_{K}}\psi_{2}|_{I_{K}}^{{-1}}=\prod_{\tau\in J}\omega_{\tau}^{n_{\tau}}\prod_{\tau\in J^{c}}\omega_{\tau}^{-n_{\tau}} $$for some $J\subset S$ and $2\leq n_{\tau}\leq p-2$. Substituting, we obtain $$\omega_{\tau_{j-1}}\prod_{\tau\in J}\omega_{\tau}^{n_{\tau}}=\omega_{\tau_{j}}^{b_{\tau_{j}}-1}\prod_{\tau\in J^{c}}\omega_{\tau}^{n_{\tau}}.$$If $\tau_{j}\in J$ then we can immediately equate coefficients of $\omega_{\tau_{j-1}}$ and deduce a contradiction. If not, then because $n_{\tau_{j}}+b_{\tau_{j}}<2p$ we see that we can still equate coefficients of $\omega_{\tau_{j-1}}$ to obtain a contradiction.

The proof in the irreducible
case is very similar, and rather simpler, as less ``carrying'' is possible. In fact, the argument gives the stronger result that $\rhobar|_{G_{F_{v}}}$ is compatible with $\sigma^{v}_{\vec{a},\vec{b}}$ for all $v$. A proof is given in the proof of Proposition 5.11 of \cite{sch06}; note that \cite{sch06} works in the setting of \cite{bdj} (using indefinite quaternion algebras), but the proof of Proposition 5.11 is purely local (using Raynaud's theory of finite flat group schemes of type $(p,\dots,p)$ in place of the Breuil module calculations used in this paper), and applies equally well in our setting.\end{proof}

The following theorem is due to Michael Schein in the case that
$\rhobar|_{G_{F_v}}$ is irreducible for all places $v|p$ (see \cite{sch06}).
\begin{thm}If $\overline{\rho}$ is modular of weight $\sigma$, and $\sigma$ is regular, then $\sigma\in
W(\overline{\rho})$.
\end{thm}
\begin{proof}Suppose that $\sigma=\otimes_v
  \sigma^v_{\vec{a},\vec{b}}$, so that we need to show that
  $\sigma^v_{\vec{a},\vec{b}}\in W_v(\overline{\rho})$ for all
  $v|p$. By Lemma \ref{compat}, $\sigma^v_{\vec{a},\vec{b}}$ is
  compatible with $\overline{\rho}|_{G_{F_v}}$, via $J$, say. If
  $\rhobar|_{G_{F_v}}$ is irreducible, we are done, so assume that it
  is reducible. By Lemma \ref{lem:local langlands version of
    lifting}, $\rhobar|_{G_{F_{v}}}$ has a lift to a potentially
  Barsotti-Tate representation of type $\tau_{\vec{a},\vec{b},J}$. By
  definition, this is, up to an unramified twist, a lift of type $J$. By Proposition
  \ref{prop:liftimpliesmodel}, $\rhobar|_{G_{F_{v}}}$ has a model of
  type $J$.  Twisting, we may without loss of generality suppose that
  each $a_\tau=0$. Then by Proposition \ref{h1f}, and the definition
  of $W_v(\overline{\rho})$, we see that
  $\sigma^v_{\vec{a},\vec{b}}\in W_v(\overline{\rho})$, as required.\end{proof}

\begin{thm}If $\sigma\in W({\rhobar})$ is a regular weight, and $\sigma$ is non-ordinary, then $\overline{\rho}$ is modular of weight $\sigma$. If $\sigma\in W({\rhobar})$ is regular, and $\sigma$ is partially ordinary of type $I$ and $\overline{\rho}$ has a partially ordinary modular lift of type $I$ then $\overline{\rho}$ is modular of weight $\sigma$.
\end{thm}
\begin{proof}Suppose that $\sigma=\otimes_v \sigma^v_{\vec{a},\vec{b}}$, so that $\sigma^v_{\vec{a},\vec{b}}\in
W_v(\overline{\rho})$ for all $v|p$. Firstly, we note that (by the
definition of $W_v(\overline{\rho})$) $\sigma^v_{\vec{a},\vec{b}}$ is compatible
with $\overline{\rho}|_{G_{F_v}}$, via $J_v$, say.

Consider firstly the case where $\overline{\rho}|_{G_{F_v}}$ is reducible. We claim that $\overline{\rho}|_{G_{F_v}}$ has a model of type $J_v$. To see this, we may twist, and without loss of generality
suppose that $a_\tau=0$ for all $\tau$, so that $\rhobar|_{G_{F_{v}}}\sim
\bigl(\begin{smallmatrix}\psi_1&*\\0&\psi_2\end{smallmatrix}\bigr)$, with
$\psi_1|_{I_{F_v}}=\prod_{\tau\in J_v}\omega_{\tau}^{b_\tau}$,
$\psi_2|_{I_{F_v}}=\prod_{\tau\notin J_v}\omega_{\tau}^{b_\tau}$. Now, by
Proposition \ref{h1f} (and the definition of $W(\overline{\rho}_v)$) $\overline{\rho}|_{G_{F_v}}$ has a model of type $J_v$, as required. Then Theorem \ref{localdef} shows that $\overline{\rho}|_{G_{F_v}}$ has a
potentially Barsotti-Tate deformation of type $\tau_{\vec{a},\vec{b},J_v}$.

If $\overline{\rho}|_{G_{F_v}}$ is irreducible, then Theorem \ref{localdef2}
shows that shows that $\overline{\rho}|_{G_{F_v}}$ has a
potentially Barsotti-Tate deformation of type $\tau'_{\vec{a},\vec{b},J_v}$.

By Corollary 3.1.7 of \cite{gee061} there is a modular lift $\rho:G_F\to
\GL_2(\Qpbar)$ of $\overline{\rho}$ which is potentially Barsotti-Tate of type
$\tau_{\vec{a},\vec{b},J_v}$ for each $v|p$ for which $\rhobar|_{G_{F_{v}}}$ is reducible, and of type
$\tau'_{\vec{a},\vec{b},J_v}$ for each $v|p$ for which $\rhobar|_{G_{F_{v}}}$ is irreducible. Then by Lemma \ref{lem:local langlands version of lifting}, $\overline{\rho}$ is modular of weight $\sigma'$ for some
Jordan-H\"{o}lder constituent $\sigma'$ of the reduction modulo $p$ of $\otimes_v
I_v$, where $I_{v}=I_{\vec{a},\vec{b},J_v}$ if $\rhobar|_{G_{F_{v}}}$ is reducible, and $I_{v}=I'_{\vec{a},\vec{b},J_v}$ otherwise. The result then
follows from Propositions \ref{typeswts} and \ref{typeswts2}, Remarks \ref{weaklyregular} and \ref{weaklyregular2}, and Lemma
\ref{compat}.\end{proof}

\bibliographystyle{amsalpha}
\bibliography{tobybib}

\def\cprime{$'$}
\providecommand{\bysame}{\leavevmode\hbox to3em{\hrulefill}\thinspace}
\providecommand{\MR}{\relax\ifhmode\unskip\space\fi MR }
\providecommand{\MRhref}[2]{%
  \href{http://www.ams.org/mathscinet-getitem?mr=#1}{#2}
}
\providecommand{\href}[2]{#2}
\begin{thebibliography}{BCDT01}

\bibitem[AS86]{as86}
Avner Ash and Glenn Stevens, \emph{Modular forms in characteristic {$l$} and
  special values of their {$L$}-functions}, Duke Math. J. \textbf{53} (1986),
  no.~3, 849--868.

\bibitem[BCDT01]{bcdt}
Christophe Breuil, Brian Conrad, Fred Diamond, and Richard Taylor, \emph{On the
  modularity of elliptic curves over {$\bold Q$}: wild 3-adic exercises}, J.
  Amer. Math. Soc. \textbf{14} (2001), no.~4, 843--939 (electronic).

\bibitem[BDJ10]{bdj}
Kevin Buzzard, Fred Diamond, and Frazer Jarvis, \emph{On {S}erre's conjecture
  for mod {$l$} {G}alois representations over totally real fields}, to appear
  in Duke Mathematical Journal, 2010.

\bibitem[BM02]{bm}
Christophe Breuil and Ariane M{\'e}zard, \emph{Multiplicit\'es modulaires et
  repr\'esentations de {${\rm GL}\sb 2({\bf Z}\sb p)$} et de {${\rm
  Gal}(\overline{\bf Q}\sb p/{\bf Q}\sb p)$} en {$l=p$}}, Duke Math. J.
  \textbf{115} (2002), no.~2, 205--310, With an appendix by Guy Henniart.

\bibitem[Bre98]{MR1621389}
Christophe Breuil, \emph{Construction de repr\'esentations {$p$}-adiques
  semi-stables}, Ann. Sci. \'Ecole Norm. Sup. (4) \textbf{31} (1998), no.~3,
  281--327. \MR{MR1621389 (99k:14034)}

\bibitem[Bre99]{MR1681105}
\bysame, \emph{Repr\'esentations semi-stables et modules fortement divisibles},
  Invent. Math. \textbf{136} (1999), no.~1, 89--122. \MR{MR1681105
  (2000c:14024)}

\bibitem[Bre00]{MR1804530}
\bysame, \emph{Groupes {$p$}-divisibles, groupes finis et modules filtr\'es},
  Ann. of Math. (2) \textbf{152} (2000), no.~2, 489--549. \MR{MR1804530
  (2001k:14087)}

\bibitem[Car09]{carusomodp}
Xavier Caruso, \emph{{$\mathbb{F}_p$}-repr\'esentations semi-stables},
  preprint.

\bibitem[CDT99]{cdt}
Brian Conrad, Fred Diamond, and Richard Taylor, \emph{Modularity of certain
  potentially {B}arsotti-{T}ate {G}alois representations}, J. Amer. Math. Soc.
  \textbf{12} (1999), no.~2, 521--567.

\bibitem[CL09]{MR2543474}
Xavier Caruso and Tong Liu, \emph{Quasi-semi-stable representations}, Bull.
  Soc. Math. France \textbf{137} (2009), no.~2, 185--223. \MR{MR2543474}

\bibitem[Dia07]{dia05}
Fred Diamond, \emph{A correspondence between representations of local {G}alois
  groups and {L}ie-type groups}, {$L$}-functions and {G}alois representations,
  London Math. Soc. Lecture Note Ser., vol. 320, Cambridge Univ. Press,
  Cambridge, 2007, pp.~187--206. \MR{MR2392355 (2009d:11085)}

\bibitem[FL82]{MR707328}
Jean-Marc Fontaine and Guy Laffaille, \emph{Construction de repr\'esentations
  {$p$}-adiques}, Ann. Sci. \'Ecole Norm. Sup. (4) \textbf{15} (1982), no.~4,
  547--608 (1983). \MR{MR707328 (85c:14028)}

\bibitem[Gee06]{gee052}
Toby Gee, \emph{A modularity lifting theorem for weight two {H}ilbert modular
  forms}, Math. Res. Lett. \textbf{13} (2006), no.~5-6, 805--811.

\bibitem[Gee07]{gee051}
\bysame, \emph{Companion forms over totally real fields. {II}}, Duke Math. J.
  \textbf{136} (2007), no.~2, 275--284.

\bibitem[Gee08]{gee061}
\bysame, \emph{Automorphic lifts of prescribed types}, preprint, 2008.

\bibitem[GS09]{geesavquatalg}
Toby Gee and David Savitt, \emph{Serre weights for quaternion algebras},
  preprint, 2009.

\bibitem[Kis08a]{kis04}
Mark Kisin, \emph{Moduli of finite flat group schemes, and modularity}, to
  appear in Annals of Mathematics (2008).

\bibitem[Kis08b]{kis06}
\bysame, \emph{Potentially semi-stable deformation rings}, J. Amer. Math. Soc.
  \textbf{21} (2008), no.~2, 513--546. \MR{MR2373358 (2009c:11194)}

\bibitem[Kis09]{kis07}
\bysame, \emph{Modularity of 2-adic {B}arsotti-{T}ate representations}, Invent.
  Math. \textbf{178} (2009), no.~3, 587--634. \MR{MR2551765}

\bibitem[Sav05]{sav04}
David Savitt, \emph{On a conjecture of {C}onrad, {D}iamond, and {T}aylor}, Duke
  Math. J. \textbf{128} (2005), no.~1, 141--197.

\bibitem[Sav08]{sav06}
\bysame, \emph{Breuil modules for {R}aynaud schemes}, J. Number Theory
  \textbf{128} (2008), no.~11, 2939--2950. \MR{MR2457845 (2010a:14083)}

\bibitem[Sch08]{sch06}
Michael Schein, \emph{Weights of {G}alois representations associated to
  {H}ilbert modular forms}, Journal Reine Angew. Math. \textbf{622} (2008),
  57--94.

\bibitem[Tay06]{taymero}
Richard Taylor, \emph{On the meromorphic continuation of degree two
  {$L$}-functions}, Doc. Math. (2006), no.~Extra Vol., 729--779 (electronic).

\end{thebibliography}

{\footnotesize \textit{E-mail address:} {\tt tgee@math.harvard.edu}}

\vskip 0.2cm

 {\footnotesize \sc Department of Mathematics, Harvard University}

\markleft{} \pagestyle{myheadings}

\end{document}